\pdfoutput=1
\documentclass[11pt]{amsart}
\usepackage[margin=1.1in]{geometry}
\usepackage{amsmath,amssymb,amsthm,mathtools}
\usepackage[expansion=false]{microtype}
\usepackage{enumitem}
\usepackage{xcolor}
\usepackage{tikz}
\usetikzlibrary{arrows.meta,calc,decorations.pathreplacing,decorations.markings,shapes.geometric}
\usepackage[colorlinks=true,linkcolor=blue!50!black,citecolor=blue!50!black,urlcolor=blue!50!black]{hyperref}
\hypersetup{pdftitle={Nash-Williams' orientation theorems for infinite graphs},
  pdfauthor={Amena Assem, Mohammad A. Mahmoud, Bobby Miraftab}}
\usepackage[capitalize,nameinlink]{cleveref}

\setlist{itemsep=1pt,topsep=3pt,parsep=0pt}
\newcommand{\para}[1]{\par\smallskip\noindent\textbf{#1}\enspace\ignorespaces}
\allowdisplaybreaks

\theoremstyle{plain}
\newtheorem{mainthm}{Theorem}
\newtheorem{theorem}{Theorem}[section]
\newtheorem{lemma}[theorem]{Lemma}
\newtheorem{proposition}[theorem]{Proposition}
\newtheorem{corollary}[theorem]{Corollary}
\theoremstyle{definition}
\newtheorem{definition}[theorem]{Definition}
\newtheorem{construction}[theorem]{Construction}
\newtheorem{example}[theorem]{Example}
\theoremstyle{remark}
\newtheorem{remark}[theorem]{Remark}

\crefname{mainthm}{Theorem}{Theorems}
\crefname{theorem}{Theorem}{Theorems}
\crefname{lemma}{Lemma}{Lemmas}
\crefname{proposition}{Proposition}{Propositions}
\crefname{corollary}{Corollary}{Corollaries}
\crefname{definition}{Definition}{Definitions}
\crefname{construction}{Construction}{Constructions}
\crefname{example}{Example}{Examples}
\crefname{remark}{Remark}{Remarks}
\crefname{figure}{Figure}{Figures}
\crefname{section}{Section}{Sections}
\crefname{subsection}{Section}{Sections}
\crefname{appendix}{Appendix}{Appendices}

\newcommand{\bd}{\delta}
\newcommand{\lam}{\lambda}
\newcommand{\dout}{d^{+}}
\newcommand{\din}{d^{-}}
\DeclareMathOperator{\Br}{Br}
\DeclareMathOperator{\depth}{depth}
\DeclareMathOperator{\tail}{tail}
\DeclareMathOperator{\head}{head}
\DeclareMathOperator{\parent}{parent}
\newcommand{\cW}{\mathcal W}
\newcommand{\cK}{\mathcal K}
\newcommand{\cU}{\mathcal U}
\newcommand{\cY}{\mathcal Y}
\newcommand{\cX}{\mathcal X}
\newcommand{\cT}{\mathcal T}
\newcommand{\cP}{\mathcal P}
\newcommand{\cC}{\mathcal C}
\newcommand{\cA}{\mathcal A}
\newcommand{\NN}{\mathbb N}
\newcommand{\Lad}{\mathbb L}
\newcommand{\Dlad}{D_{\mathbb L}}

\tikzset{
  vx/.style={circle,fill=black,inner sep=1.25pt},
  tm/.style={rectangle,fill=black,inner sep=1.9pt},
  gd/.style={diamond,draw=black,fill=white,thick,inner sep=1.6pt},
  ar/.style={-{Stealth[length=2mm]},semithick},
  arthick/.style={-{Stealth[length=2.2mm]},line width=1.4pt},
  cyc/.style={-{Stealth[length=2.2mm]},line width=1.4pt,red!75!black},
  cyccore/.style={-{Stealth[length=2.2mm]},line width=0.9pt,red!75!black,dash pattern=on 2.6pt off 1.8pt},
  core/.style={fill=gray!22,draw=none},
  mid/.style={semithick,postaction={decorate},decoration={markings,mark=at position 0.56 with {\arrow{Stealth[length=2mm]}}}},
  midthick/.style={line width=1.4pt,postaction={decorate},decoration={markings,mark=at position 0.56 with {\arrow{Stealth[length=2.3mm]}}}},
  midcyc/.style={line width=1.4pt,red!75!black,postaction={decorate},decoration={markings,mark=at position 0.56 with {\arrow{Stealth[length=2.3mm]}}}},
}

\title{Nash-Williams' orientation theorems for infinite graphs}
\author{Amena Assem}
\address[A.~Assem]{York University, Canada}
\author{Mohammad A.~Mahmoud}
\address[M.~A.~Mahmoud]{University of Toronto, Canada}
\author{Bobby Miraftab}
\address[B.~Miraftab]{Carleton University, Canada}
\date{}

\begin{document}

\begin{abstract}
Nash-Williams proved in 1960 that every finite multigraph has a \emph{well-balanced} orientation: for all vertices
$x$ and $y$, the maximum number of arc-disjoint directed paths from $x$ to $y$ is at least half, rounded down, of the
maximum number of edge-disjoint $x$--$y$ paths. In particular every finite $2k$-edge-connected multigraph has a
$k$-arc-connected orientation. He announced an extension to infinite graphs whose details, as far as is known, never
appeared, and the
two statements for infinite graphs have since been known as the strong and the weak orientation conjecture of
Nash-Williams. We prove both, for multigraphs of arbitrary cardinality.

For multigraphs with countably many edges we obtain more: a well-balanced orientation that is smooth (in-degree and
out-degree differ by at most one at every vertex of finite degree) and, if the multigraph is $2$-edge-connected,
has infinitely many arcs in each direction across every infinite cut. Such an orientation is obtained by
K\H{o}nig's lemma from a sequence of finite orientation problems, each solved by the theorem of Kir\'aly and Szigeti
on simultaneous well-balanced orientations. The infinite cuts are controlled by a canonical tree of
$2$-edge-connected pieces with finite cores: once the cores are strongly connected and every finite cut is crossed in both
directions, every arc lies on a directed cycle that uses only edges owned by one piece and edges of the cores of its children, and a counting argument of Thomassen does the rest.
Neither local finiteness nor ends play any role.
\end{abstract}

\maketitle

\section{Introduction}\label{sec:intro}

\subsection{Orientations of finite graphs}\label{sec:intro-finite}

An \emph{orientation} of a graph replaces each edge by an arc between the same two vertices. In 1939 Robbins
\cite{Rob39}, motivated by a problem of traffic control, proved that a finite graph with at least two vertices has a
strongly connected orientation if and only if it is $2$-edge-connected. Nash-Williams \cite{NW60} found the generalization to higher
connectivity: \emph{every finite $2k$-edge-connected multigraph has a $k$-arc-connected orientation}, that is, an
orientation with $k$ pairwise arc-disjoint directed paths from any vertex to any other. Since every cut of a
$k$-arc-connected digraph contains at least $k$ arcs in each direction, edge-connectivity $2k$ is also necessary.
(Lov\'asz's and Mader's splitting-off theorems give other proofs of this statement; see \cite[p.~684]{KSz06} and
\cite[p.~602]{Th16}.)

Nash-Williams derived this from a local statement. For distinct vertices $x,y$ of a graph $G$ let $\lam_G(x,y)$ be
the maximum number of pairwise edge-disjoint $x$--$y$ paths, and for an orientation $D$ of $G$ let $\lam_D(x,y)$ be
the maximum number of pairwise arc-disjoint directed paths from $x$ to $y$. The orientation $D$ is
\emph{well-balanced} if
\[\lam_D(x,y)\ \ge\ \lfloor\lam_G(x,y)/2\rfloor\qquad\text{for all }x\ne y.\]
Nash-Williams' Theorem~1 in \cite{NW60} states that \emph{every finite graph has a well-balanced orientation} (he
called such orientations \emph{admissible}). If $G$ is $2k$-edge-connected, then $\lam_G(x,y)\ge 2k$ for all $x\ne y$
by Menger's theorem, so a well-balanced orientation is $k$-arc-connected. His proof rests on the theorem on
\emph{odd-vertex pairings} \cite[Thm.~2]{NW60}: the vertices of odd degree can be joined in pairs by a set $P$ of
new edges in such a way that $|\bd_G(S)|-|\bd_P(S)|\ge2\lfloor\lam_G(x,y)/2\rfloor$ for all distinct vertices
$x,y$ and every vertex set $S$ containing $x$ but not $y$. Every Eulerian orientation of the graph $G+P$, restricted to $G$, is then well-balanced
\cite[\S3]{NW60} (see also the proof of \cite[Thm.~5.1]{AKP25}), and it is also \emph{smooth}: in-degree and
out-degree differ by at most one at every vertex. A smooth well-balanced orientation is called
\emph{best-balanced}, so every finite graph has a best-balanced orientation \cite[Thm.~3]{KSz06}, where this is
credited to \cite{NW60,NW69}.

\subsection{The infinite problem}\label{sec:intro-infinite}

In infinite graphs, cuts and local connectivities may be infinite. We read $\lam_G(x,y)$ and $\lam_D(x,y)$ as the
suprema of the cardinalities of families of pairwise edge-disjoint $x$--$y$ paths, respectively of pairwise
arc-disjoint directed $x\to y$ paths (paths are always finite); these suprema are attained and equal the least size
of a separating cut (\cref{lem:lam}). Following Pitz and Stegemann \cite[Def.~2.1]{PS26}, an orientation $D$ of a graph
$G$ is \emph{well-balanced} if for all distinct vertices $x,y$
\begin{equation}\label{eq:wb}
\lam_D(x,y)\ \ge\ \begin{cases}\lfloor\lam_G(x,y)/2\rfloor&\text{if }\lam_G(x,y)<\infty,\\[1pt] \lam_G(x,y)&\text{otherwise,}\end{cases}
\end{equation}
and $D$ is \emph{smooth} if $|\dout_D(v)-\din_D(v)|\le1$ at every vertex $v$ of finite degree.

In the paper that proved the finite theorems, Nash-Williams wrote: ``An extension of Theorem~1 to infinite graphs has
been obtained, but the details, being somewhat heavy, are deferred to a possible future paper'' \cite[p.~556]{NW60}.
As far as is known, no such paper appeared \cite[p.~602]{Th16}. The infinite versions of his two theorems became known as the \emph{weak} and the
\emph{strong Nash-Williams orientation conjecture}: every $2k$-edge-connected graph has a $k$-arc-connected
orientation, and every graph has a well-balanced orientation. This paper proves both (\cref{thm:weak,thm:strong}
below).

\subsection{What changes in infinite graphs}\label{sec:intro-hard}

Many finite results extend to infinite graphs by compactness. For orientations the natural compact space is
$\{0,1\}^{E(G)}$, one coordinate per edge, and a condition that depends on the directions of finitely many edges
defines a closed (indeed clopen) set. That a \emph{finite} cut contains at least $k$ arcs in each direction is such
a condition. Thomassen remarked about two related orientation problems for infinite graphs that they ``have
affirmative answers if we focus only on finite cuts. It is the infinite cuts that are problematic''
\cite[p.~618]{Th16}.

\begin{figure}[t]
\centering
\begin{tikzpicture}[x=1.12cm,y=1.25cm]
\foreach \n in {0,...,7}{
  \node[vx,label={[font=\scriptsize,yshift=-1pt]above:$a_{\n}$}] (a\n) at (\n,1) {};
  \node[vx,label={[font=\scriptsize,yshift=1pt]below:$b_{\n}$}] (b\n) at (\n,0) {};
}
\foreach \n [count=\m from 1] in {0,...,6}{
  \draw[ar] (a\n) -- (a\m);
  \draw[ar] (b\m) -- (b\n);
}
\foreach \n in {1,...,7}{ \draw[ar] (a\n) -- (b\n); }
\draw[arthick,red!75!black] (b0) -- (a0);
\node at (7.75,1) {$\cdots$};
\node at (7.75,0) {$\cdots$};
\draw[dashed,gray] (-0.45,0.5) -- (8.05,0.5);
\node[font=\scriptsize,anchor=west] at (8.1,0.5) {$\bd(\{a_n:n\ge0\})$};
\end{tikzpicture}
\caption{The one-way ladder $\Dlad$: the top rail is directed to the right, the bottom rail to the left, every rung
downwards except the first. It is strongly connected, but the infinite cut $\bd(\{a_n:n\ge0\})$ (dashed), which
consists of all rungs, is crossed only once from bottom to top (thick).}\label{fig:ladder}
\end{figure}
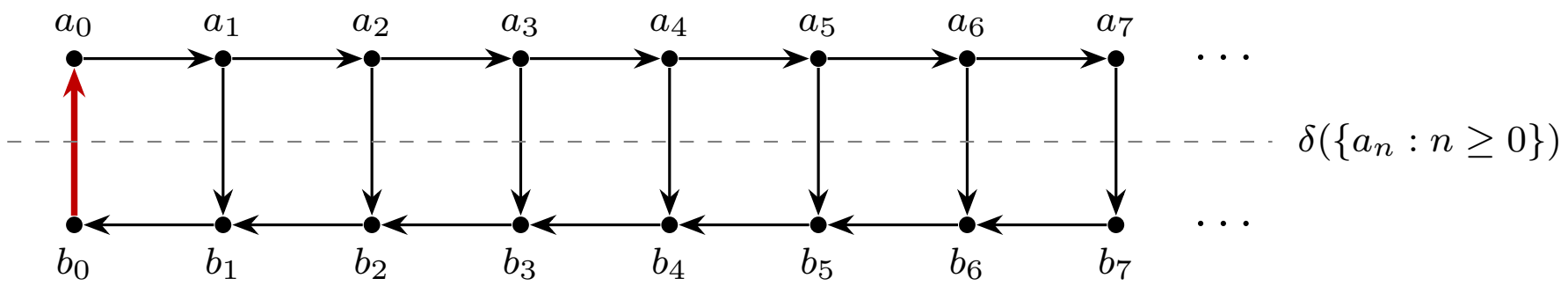

Infinite cuts are not finite conditions. If $\bd(A)$ is an infinite cut, the set of orientations with an arc
entering $A$ is open but not closed: in a limit, the entering arcs may escape to infinity.
\Cref{fig:ladder} shows the ladder $\Lad$ with an orientation $\Dlad$ that is strongly connected, but crosses the
infinite cut formed by all rungs only once from the bottom to the top. It is a limit of orientations each of which
has infinitely many arcs in each direction across every infinite cut (\cref{rem:notclosed}). On the ladder with
$k$ parallel copies of every edge, $k$ copies of $\Dlad$ form a $k$-arc-connected orientation in which this cut is
crossed only $k$ times upwards, and no reorientation that changes the out-degrees at only finitely many
vertices---for instance, the reversal of any family of arc-disjoint directed cycles---can produce infinitely many upward
rungs (\cref{prop:fibre}). Infinite cuts must
therefore be controlled while the orientation is being built.

Two other finite tools are unavailable as well. In a finite $4k$-edge-connected graph the tree-packing theorem of
Nash-Williams and Tutte provides $2k$ edge-disjoint spanning trees, and orienting $k$ of them towards a common root and
$k$ away from it gives a $k$-arc-connected orientation \cite[\S1]{AKP25}. But Aharoni and Thomassen \cite{AT89}
constructed countable graphs of arbitrarily large finite edge-connectivity that do not even have two edge-disjoint
spanning trees. And the obvious way to reduce to finite graphs---contract the infinite parts that remain after
removing a finite set---loses information about how these parts are connected at infinity, that is, about the ends
of the graph. How this has been dealt with is described in \cref{sec:history}, after the statement of our results.

\subsection{Results}\label{sec:results}

We call an orientation $D$ of a graph $G$ \emph{two-way infinite} if every infinite cut of $G$ contains infinitely
many arcs of $D$ in each direction. All graphs are multigraphs without loops, and parallel classes, degrees and vertex
sets may have any cardinality.

\begin{mainthm}[weak orientation theorem]\label{thm:weak}
Let $k\ge1$ be an integer. Every $2k$-edge-connected multigraph has a $k$-arc-connected orientation.
\end{mainthm}

\begin{mainthm}[strong orientation theorem]\label{thm:strong}
Every multigraph has a well-balanced orientation.
\end{mainthm}

\begin{mainthm}[countably many edges]\label{thm:ctble}
Let $G$ be a multigraph with countably many edges.
\begin{enumerate}[label=\textup{(\roman*)}]
\item $G$ has a smooth well-balanced orientation.
\item If $G$ is $2$-edge-connected, then $G$ has a smooth well-balanced orientation $D$ that is two-way infinite.
Such a $D$ is strongly connected, satisfies $\dout_D(v)=\din_D(v)=\aleph_0$ at every vertex $v$ of infinite degree,
and is $k$-arc-connected whenever $G$ is $2k$-edge-connected.
\end{enumerate}
\end{mainthm}

\Cref{thm:weak,thm:strong} are the weak and the strong orientation conjecture of Nash-Williams. \Cref{thm:strong}
implies \cref{thm:weak} (\cref{lem:wb-to-k}), but we also prove \cref{thm:weak} by Thomassen's Zorn's lemma argument
\cite[\S8]{Th16}, which avoids the bond-faithful decompositions behind the reduction of Pitz and Stegemann
\cite[Thm.~1.2]{PS26}. \Cref{thm:ctble} is the heart of the paper: all the new work lies in the countable case, and \cref{thm:weak,thm:strong} follow from it by the two known reductions (\cref{sec:all}). Its part~(ii) strengthens Thomassen's Theorem~11
\cite{Th16} (see \cref{sec:history}) for graphs with countably many edges.
No hypothesis on local finiteness or on the ends is needed anywhere, and multigraphs with infinite parallel classes are
allowed throughout. Smoothness and the two-way infinite property are not claimed beyond countably many edges
(\cref{rem:smooth-open}).

\subsection{Previous work}\label{sec:history}

We describe the contributions to the problem in chronological order, together with the obstacles they met.

\para{Early results.} The case $k=1$ of the weak conjecture is Robbins' theorem for finite graphs
\cite{Rob39}; Egyed \cite{Egy41} extended it to all infinite graphs in 1941 (see also \cite[p.~556]{NW60};
\cite{AKP25} describes Egyed's proof as a Zorn's lemma argument). Nash-Williams' announcement of 1960 quoted above has
been followed by two accounts. According to \cite{PS26,AKP25}, he retracted the claim in \cite[\S8]{NW69}. According
to Thomassen \cite[p.~602]{Th16}, following Bar\'at and Kriesell \cite{BK10}, the infinite statement was stated by
Nash-Williams in his 1967 survey \cite{NW67} and apparently never published. In \cite{NW69} Nash-Williams also
formulated a \emph{simultaneous} orientation theorem: every subgraph $H$ of a finite graph $G$ has a
best-balanced orientation that extends to a best-balanced orientation of $G$ (see \cite[p.~685 and
Thm.~4]{KSz06}). In 1989 Aharoni and Thomassen \cite{AT89} showed that tree packing fails in infinite graphs, as
described above. Thomassen conjectured that \emph{some} function $f(k)$ of the edge-connectivity forces a
$k$-arc-connected orientation (see \cite[p.~602]{Th16}).

\para{Simultaneous well-balanced orientations.} Kir\'aly and Szigeti \cite[Thm.~8]{KSz06} proved in 2006 that
for every partition of the edge set of a finite graph there is a best-balanced orientation of the whole graph that is
best-balanced on every part as well; this generalizes the theorem of \cite{NW69} just mentioned, and its proof is
a short argument from Nash-Williams' odd-vertex-pairing theorem. Bern\'ath, Iwata, Kir\'aly, Kir\'aly and Szigeti
\cite{BIKKS08} discuss such simultaneous statements and their limits; for nested subgraphs they fail
\cite[Counter-example~6]{BIKKS08}. This theorem is the finite engine of our proof.

\para{Thomassen's theorem.} In 2016 Thomassen \cite{Th16} proved that $8k$-edge-connectivity suffices
\cite[Thm.~7]{Th16} (parallel edges of finite multiplicity are allowed there), the first bound of any kind. His proof
set up the framework of the subsequent work on the weak conjecture \cite{As25,AKP25,AR25}. A vertex set $B$ is \emph{boundary-linked} if the boundary edges of $B$ are
the first edges of pairwise edge-disjoint rays in $G[B]\cup\bd(B)$ that belong to a common end. Thomassen showed
that in a connected locally finite graph, for every vertex set $A_0$ with finite boundary, $V(G)\setminus A_0$ splits
into finitely many singletons and boundary-linked sets with finite boundaries \cite[Thm.~1]{Th16}. Contracting these sets and \emph{lifting} (splitting
off) pairs of edges at the contracted vertices, with the help of his theorem that the ``bad graph'' (the complement
of the lifting graph) is disconnected in the Eulerian case \cite[Thm.~2]{Th16}, he found immersions of finite $2k$-edge-connected Eulerian graphs in $4k$-edge-connected
graphs \cite[Thm.~4]{Th16}, and he oriented an increasing sequence of such finite immersions. Two reductions extend
the result from locally finite to countable graphs, by splitting vertices \cite[\S7, Thm.~9]{Th16}, and from
countable to arbitrary graphs, by Zorn's lemma with \cite[Thm.~10]{Th16} supplying countable $m$-edge-connected
subgraphs \cite[\S8]{Th16}. The same paper contains two observations that are relevant here. First, Theorem~11 there:
\emph{every infinite $2$-edge-connected graph has a strongly connected orientation in which every infinite cut has
infinitely many arcs in both directions}. Its proof assigns to every edge a directed cycle and counts how often each
edge is ``demanded'' by these cycles; it is given in detail for locally finite graphs, extended to countable graphs
via Theorem~9 there, and the uncountable case is left to the reader \cite[p.~619]{Th16}. Second, by Nash-Williams'
cycle-decomposition theorem, a countable $2k$-edge-connected graph without finite odd cuts has a $k$-arc-connected orientation:
orient each cycle of a decomposition into finite cycles cyclically \cite[p.~601]{Th16}.

\para{From $8k$ to $2k$: lifting, and countably many ends.} In 2020 the preprint \cite{KP20} claimed the weak
conjecture for all infinite multigraphs; this claim was later retracted. In 2022 the \emph{lifting graph} $L(G,s,k)$, the
graph on the edges at a vertex $s$ in which two edges are adjacent if lifting them preserves local edge-connectivity $k$
between the other vertices, was analysed in \cite{As22}. Its structure was determined there for all parities of $\deg(s)$ and $k$, building on results of
Ok, Richter and Thomassen \cite{ORT16}, and extending to odd $\deg(s)$ results that Jord\'an \cite{Jo03} had found
independently; in particular, for $k\ge2$ and $(s,k)$-edge-connected $G$ with $\deg(s)>4$, if $L(G,s,k)$ is
disconnected then $\deg(s)$ is odd and $L(G,s,k)$ consists of an isolated vertex and a balanced complete bipartite
graph \cite[Thm.~1.5(i)]{As22}. The same paper related its maximal independent sets to the dangerous sets of $G$. With these
results Thomassen's bound was lowered to $4k$ for one-ended locally finite graphs \cite[Thm.~6.1]{As25}; the new ingredient is a lifting
sequence at a contracted vertex of odd degree, including the exceptional configuration just described
\cite[Thm.~2.3 and Lemma~3.1]{As25}. The optimal bound $2k$ for locally finite graphs with countably many ends was then
proved in \cite[Thm.~6.1]{AKP25}, which builds on work from \cite{KP20} together with the new results on the
lifting graph of \cite{As22}. A key step is a decomposition in which the boundary-linked sets are exactly
the components of $G-A$ for a finite set $A$ \cite[Thm.~2.2]{AKP25}. It needs countably many ends: the \emph{twin
binary tree} \cite[Fig.~1]{AKP25}, a locally finite graph with uncountably many ends, shows that this hypothesis cannot
be dropped. Further ingredients are an immersion theorem with auxiliary vertices of degree three \cite[Thm.~3.2]{AKP25}
and the extension of consistently oriented open or closed Eulerian subgraphs to well-balanced orientations of finite
graphs \cite[Thm.~5.1]{AKP25}; substituting the latter into Thomassen's proof gives the bound $4k$ for all graphs (the
remark opening \cite[\S6]{AKP25}, credited to Pitz in \cite[\S5]{AR25}).

\para{The strong conjecture.} Pitz and Stegemann \cite{PS26} (2024) proved the strong conjecture for
rayless graphs \cite[Thm.~1.1]{PS26} and reduced it to countable graphs \cite[Thm.~1.2]{PS26}, using Laviolette's
theorem (2005) that every graph has a bond-faithful decomposition into countable connected subgraphs \cite{La05}, in
the form reproved by Thomassen in 2017 without Nash-Williams' cycle-decomposition theorem \cite[Thm.~3]{Th17}.

\para{$2k+1$ for all graphs.} In 2025 it was proved in \cite{AR25} that Huck's theorem \cite{Hu91} holds for all
infinite graphs---for odd $k$, every $(k+1)$-edge-connected infinite graph is weakly $k$-linked \cite[Thm.~1.4]{AR25}---using a
\emph{linking fan} \cite[Prop.~3.1]{AR25} and lift sequences that are compatible with the witnessing rays of contracted
boundary-linked sets \cite[Lemma~2.5]{AR25}. The same tools give immersions of finite $2k$-edge-connected graphs in
$(2k+1)$-edge-connected locally finite graphs, with no restriction on the ends and with edges between the
boundary-linked sets allowed \cite[Thm.~5.1, Cor.~5.2]{AR25}, and hence that $2k+1$ suffices for all infinite graphs
\cite[Thm.~5.3]{AR25}. In that immersion the boundary-linked sets of odd boundary survive as vertices of degree $2k+1$
\cite[Thm.~5.1(ii)]{AR25}. In our reading, this is where the extra unit of connectivity is spent: at connectivity
exactly $2k$, a set of odd boundary may lead to the exceptional configuration of the lifting graph described above.

\para{Edge-ends.} Aurichi, Magalh\~aes and Pinto \cite{AMP26} (2025) proved the weak conjecture for graphs with countably many
\emph{edge-ends} (classes of rays under edge-disjoint linkage, which for locally finite graphs are the ends), without
assuming local finiteness \cite[Thm.~3.2]{AMP26}, and reduced the strong conjecture to locally finite graphs
\cite[Thm.~4.1]{AMP26}.

\subsection{The idea of the proof}\label{sec:idea}

We now explain the proof in plain words; the details occupy \cref{sec:roadmap,sec:tree,sec:cycles,sec:engine,sec:compact}.

\para{The demand count.} We do not try to control infinite cuts directly. Instead we use the counting argument
from Thomassen's proof of his Theorem~11. Suppose that every arc $e$ of an orientation $D$ lies on a finite directed
cycle $C_e$, and that no arc lies on infinitely many of the chosen cycles. If an infinite cut $\bd(A)$ had only
finitely many arcs entering $A$, then each of the infinitely many arcs leaving $A$ would lie on a cycle that returns
to $A$ through one of these finitely many arcs---impossible, since each of them lies on only finitely many of the
cycles. So $D$ is two-way infinite (\cref{thm:dp}). Our task is thus to find, by finite means, an orientation in which
every arc lies on a directed cycle and the cycles can be chosen so that no arc is used infinitely often. Strong
connectivity is not enough: the return path of an arc may wander arbitrarily far.

\begin{figure}[t]
\centering
\begin{tikzpicture}[x=1.12cm,y=1.2cm]
\foreach \j in {0,2,4,6}{
  \fill[core,rounded corners=4pt] (\j-0.24,-0.24) rectangle (\j+1.24,1.24);
}
\foreach \n in {0,1,3,5,7}{
  \node[vx] (a\n) at (\n,1) {};
  \node[vx] (b\n) at (\n,0) {};
}
\foreach \n in {2,4,6}{
  \node[tm] (a\n) at (\n,1) {};
  \node[tm] (b\n) at (\n,0) {};
}
\foreach \n in {0,...,7}{
  \node[font=\scriptsize,above=3pt] at (a\n) {$a_{\n}$};
  \node[font=\scriptsize,below=3pt] at (b\n) {$b_{\n}$};
}
\draw[ar] (a1) -- (a0); \draw[ar] (a0) -- (b0); \draw[ar] (b0) -- (b1);
\foreach \j/\k in {2/3,4/5,6/7}{
  \draw[ar] (b\j) -- (b\k); \draw[ar] (b\k) -- (a\k); \draw[ar] (a\k) -- (a\j);
}
\draw[ar] (a4) -- (b4); \draw[ar] (a6) -- (b6);
\draw[arthick] (a3) -- (a4); \draw[arthick] (b4) -- (b3);
\draw[arthick] (a5) -- (a6); \draw[arthick] (b6) -- (b5);
\draw[cyc] (a1) -- (a2); \draw[cyccore] (a2) -- (b2); \draw[cyc] (b2) -- (b1); \draw[cyccore] (b1) -- (a1);
\draw[red!75!black,-{Stealth[length=1.6mm]},thick] (1.5,0.72) arc[start angle=90,end angle=-200,radius=0.22];
\node at (7.75,1) {$\cdots$}; \node at (7.75,0) {$\cdots$};
\begin{scope}[yshift=-0.55cm]
\node[font=\small] (R) at (0.5,-0.9) {$\Lad$};
\node[font=\small] (P2) at (2.5,-1.7) {$\Lad_{\ge2}$};
\node[font=\small] (P4) at (4.5,-2.5) {$\Lad_{\ge4}$};
\node[font=\small] (P6) at (6.5,-3.3) {$\Lad_{\ge6}$};
\node (P8) at (7.6,-3.75) {$\ddots$};
\foreach \P/\Q/\x/\y/\i in {R/P2/0.5/-1.7/0,P2/P4/2.5/-2.5/2,P4/P6/4.5/-3.3/4}{
  \draw (\P) -- (\Q);
  \pgfmathtruncatemacro{\ii}{\i+1}
  \foreach \dx/\lab in {-0.75/a_{\i},-0.3/b_{\i},0.15/a_{\ii},0.6/b_{\ii}}{
    \node[font=\tiny,inner sep=1pt] (L) at (\x+\dx-0.25,\y) {$\lab$};
    \draw[gray] (\P) -- (L);
  }
}
\draw (P6) -- (P8);
\end{scope}
\end{tikzpicture}
\caption{The ladder $\Lad$ and its piece tree. The cores (shaded) are the squares on levels $\{2j,2j+1\}$; the new
bridges (thick) are the rails between levels $2j+1$ and $2j+2$. The root $\Lad$ has five children: the single
vertices $a_0,b_0,a_1,b_1$ and the ladder $\Lad_{\ge2}$ on the levels $\ge2$, whose terminals are $a_2,b_2$ (squares);
and so on. The orientation shown has strongly connected cores and crosses every finite cut in both directions. The
cycle of the bridge $(a_1,a_2)$ (red) uses both new bridges of the root (thick) and returns through one arc of the core
of the child $\Lad_{\ge2}$ and one arc of the root core (dashed).}
\label{fig:laddertree}
\end{figure}
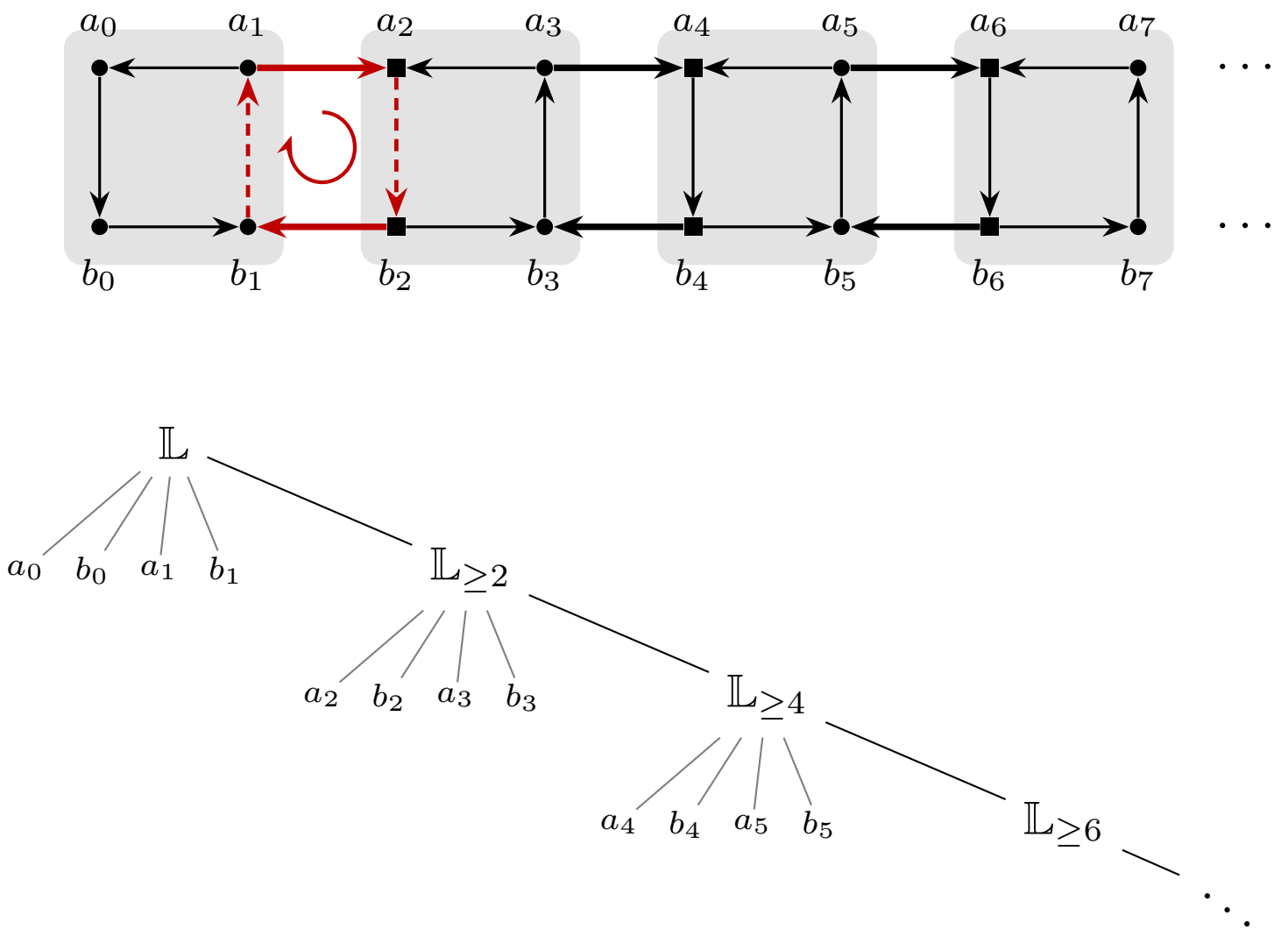

\para{The piece tree.} We confine the return paths with a canonical tree of $2$-edge-connected subgraphs, the
\emph{pieces} (\cref{sec:tree}). The root piece is $G$. Each piece $Z$ with an edge receives a finite
$2$-edge-connected \emph{core} $K(Z)\subseteq Z$ that contains the first edge of $Z$ (in a fixed enumeration of
$E(G)$) and the finitely many \emph{terminals} of $Z$, the vertices where $Z$ is attached to edges owned by its
ancestors. Deleting the edges of the core creates finitely many new bridges, the set $B(Z)$; the 2EC components
(possibly single vertices) of what remains after deleting the new bridges are the \emph{children} of $Z$, and their terminals are the vertices where they meet
the core or the new bridges. \Cref{fig:laddertree} shows the ladder. Every edge is \emph{owned}, as a core edge or a
new bridge, by exactly one piece, every piece has finitely many children, and every finite cut of $G$ is eventually a
union of pieces.

\para{Local cycles.} Suppose that an orientation $D$ crosses every finite cut in both directions---we call this
property (F)---and makes every core strongly connected---property (C). Then every arc owned by a piece $Z$ lies on a directed cycle inside the edges
owned by $Z$ and the cores of the children of $Z$ (\cref{lem:cyc}): an arc of the core returns inside the core; for
a new bridge, the children of $Z$ turn out to form a strongly connected quotient---this is where (F) enters---and
inside each child the terminals are joined by directed paths in its strongly connected core (red in
\cref{fig:laddertree}). An arc is then demanded only by arcs owned by its own piece or by the parent's new bridges,
a finite set, and the demand count shows that $D$ is two-way infinite (\cref{thm:demand}). This two-way infinite part
of the argument does not depend on $k$ or on well-balancedness: it uses only (F) and (C).

\para{The finite engine and compactness.} The properties (F) and (C), well-balancedness on finite cuts, and
smoothness at the vertices of finite degree are all conjunctions of requirements on finitely many edges each. We impose them level by level: at level
$n$ we orient the finite set $E_n$ of edges owned by pieces of depth at most $n$ on an auxiliary finite graph $H_n$, in
which every piece of depth $n+1$ is replaced by a \emph{gadget}, a new vertex joined to the terminals of the piece by
many parallel edges (\cref{fig:gadget}). The multiplicity is chosen so large that no local edge-connectivity that a
finite cut of $G$ can see is lowered. On $H_n$ we need a best-balanced orientation in which the cores of depth at most
$n$ are strongly connected; this is exactly the theorem of Kir\'aly and Szigeti, applied to the partition of $E(H_n)$
into the cores and the remaining edges (\cref{sec:engine}). One application of K\H{o}nig's lemma then yields an
orientation of $G$ that is well-balanced on finite cuts and smooth and has strongly connected cores. By the local
cycles it is two-way infinite, and the two properties together make it well-balanced (\cref{lem:cuts}). Parity, the
source of the difficulties at odd boundaries in the lifting approach, never appears outside the finite theorem.

\para{Uncountable graphs.} Graphs with uncountably many edges are handled by two known reductions: Thomassen's
Zorn's lemma argument \cite[\S8]{Th16} for the weak theorem, and the reduction theorem of Pitz and Stegemann
\cite[Thm.~1.2]{PS26} for the strong one (\cref{sec:all}).

\subsection{Organization}\label{sec:org}
\Cref{sec:prelim} fixes conventions and recalls Menger's theorem for infinite graphs (proofs in
\cref{app:menger}). \Cref{sec:roadmap} isolates the two properties---well-balanced on finite cuts, and a directed cycle
cover of finite demand---and shows that together they give a well-balanced orientation. \Cref{sec:tree} builds the
piece tree, \cref{sec:cycles} proves that strongly connected cores yield local cycles of finite demand, \cref{sec:engine}
states the finite theorem we use, and \cref{sec:compact} carries out the compactness argument over the levels and proves
\cref{thm:ctble}. \Cref{sec:all} treats arbitrary cardinalities and proves \cref{thm:weak,thm:strong}. \Cref{sec:examples}
collects examples and remarks.

\section{Preliminaries}\label{sec:prelim}

\subsection{Conventions}\label{sec:conv}

All graphs are multigraphs without loops; vertex sets, parallel classes and degrees may have any cardinality. For a
graph $G=(V,E)$ and $A\subseteq V$ we write $\bd_G(A)$ for the set of edges with exactly one end in $A$. A \emph{cut} is
a set $\bd_G(A)$ with $\emptyset\ne A\subsetneq V$; note that $\bd_G(A)=\bd_G(V\setminus A)$. For an orientation $D$ of
$G$, $\bd^+_D(A)$ and $\bd^-_D(A)$ are the sets of arcs leaving and entering $A$, and $\dout_D(v)$, $\din_D(v)$ are the
out-degree and the in-degree of a vertex $v$. For $m\in\NN$, $G$ is \emph{$m$-edge-connected} if $|V|\ge2$ and every cut
has at least $m$ edges; for $m\ge1$ this implies that $G$ is connected, since a disconnected graph has an empty cut. A
digraph $D$ is \emph{$k$-arc-connected} ($k\in\NN$) if $|V(D)|\ge2$ and $D$ has $k$ pairwise arc-disjoint directed
$x\to y$ paths for all $x\ne y$.

For $F\subseteq E$, $G-F$ is the graph $(V,E\setminus F)$, and $G[A]$ is the subgraph induced by $A\subseteq V$. A
\emph{subgraph} is a pair $(V',E')$ with $V'\subseteq V$, $E'\subseteq E$ and all ends of edges of $E'$ in $V'$. For a
subgraph $H$ and an orientation $D$ of $G$, $D[H]$ is the orientation of $H$ induced by $D$, and for $F\subseteq E$,
$D|_F$ is the restriction of $D$ to $F$. Paths and cycles are finite and have no repeated vertices; a cycle may consist
of two parallel edges. Directed paths and cycles use each of their arcs in the forward direction.

A \emph{bridge} of $G$ is an edge $b=uv$ such that $u$ and $v$ lie in different components of $G-b$; equivalently
(join a $u$--$v$ path of $G-b$ to $b$), $b$ lies on no cycle. We write $\Br(G)$ for the set of bridges, and call $G$
\emph{2EC} if it is connected and has no bridge; a single vertex is 2EC. For a connected graph with at least two
vertices, having no bridge is the same as every cut having at least two edges: if $\bd_G(A)=\{b\}$ then $b$ is a bridge,
and if $b=uv$ is a bridge then the vertex set $A$ of the component of $u$ in $G-b$ has $\bd_G(A)=\{b\}$. So a 2EC graph
with at least two vertices is the same as a $2$-edge-connected graph.

\subsection{Menger's theorem and local connectivity}\label{sec:menger}

We need Menger's theorem for finitely many paths in graphs of arbitrary cardinality. The usual augmenting path argument
proves it without any compactness; for completeness we give the proofs of the two lemmas of this subsection in
\cref{app:menger}. For a digraph $D$ and vertices $x\ne y$ let $\cA_D(x,y)$ be the family of sets $A\subseteq V(D)$
with $x\in A$ and $y\notin A$.

\begin{lemma}[Menger's theorem at a finite parameter]\label{lem:menger}
Let $k\in\NN$ and let $x\ne y$ be vertices.
\begin{enumerate}[label=\textup{(\alph*)}]
\item A digraph $D$ of any cardinality has $k$ pairwise arc-disjoint directed $x\to y$ paths if and only if
$|\bd^+_D(A)|\ge k$ for every $A\in\cA_D(x,y)$.
\item If a graph $G$ of any cardinality satisfies $|\bd_G(A)|\ge k$ for every $A\subseteq V(G)$ with $x\in A$ and
$y\notin A$, then $G$ has $k$ pairwise edge-disjoint $x$--$y$ paths.
\end{enumerate}
\end{lemma}

The local connectivities $\lam_G(x,y)$ and $\lam_D(x,y)$ were defined in \cref{sec:intro-infinite} as suprema. They are
in fact maxima, and equal to the least size of a separating cut. For graphs with countably many edges this needs no
choice principle; at arbitrary cardinality we use the Axiom of Choice.

\begin{lemma}[local connectivity is a minimum cut]\label{lem:lam}
Let $D$ be a digraph and $x\ne y$ vertices of $D$.
\begin{enumerate}[label=\textup{(\alph*)}]
\item If $D$ has countably many arcs, put $\kappa:=\min\{|\bd^+_D(A)|:A\in\cA_D(x,y)\}\in\NN\cup\{\aleph_0\}$. Then $D$
has a family of exactly $\kappa$ pairwise arc-disjoint directed $x\to y$ paths and no larger family; thus
$\lam_D(x,y)=\kappa$ and the supremum defining $\lam_D(x,y)$ is attained.
\item \textup{(Axiom of Choice.)} If $D$ has any cardinality, the set $\{|\bd^+_D(A)|:A\in\cA_D(x,y)\}$ of cardinals has
a least element $\kappa$, and the conclusion of \textup{(a)} holds.
\end{enumerate}
The same holds for a graph $G$ \textup{(}with countably many edges in \textup{(a)}, of any cardinality in
\textup{(b))}, with $\bd_G(A)$ \textup{(}$x\in A$, $y\notin A$\textup{)} in place of $\bd^+_D(A)$ and families of
edge-disjoint $x$--$y$ paths in place of families of arc-disjoint directed paths.
\end{lemma}

Consequently, at every cardinality, $\lam_G(x,y)$ and $\lam_D(x,y)$ are maxima and equal the least cardinality of a
separating cut, and condition~\eqref{eq:wb} means the same whether $\lam$ is read as a supremum, as a maximum or as a
minimum cut. In particular, an orientation is well-balanced in the sense of \eqref{eq:wb} if and only if it is
well-balanced in the sense of \cite[Def.~2.1]{PS26}, whatever the cardinality of the graph. Part~(b) of \cref{lem:lam}
is used only for this identification, in \cref{sec:PS}; the proofs of \cref{thm:weak,thm:ctble} do not need it.

By \cref{lem:menger}(a), applied to every ordered pair of distinct vertices, a digraph $D$ with $|V(D)|\ge2$ is
$k$-arc-connected if and only if
\begin{equation}\label{eq:cutform}
|\bd^+_D(A)|\ \ge\ k\qquad\text{for every }\emptyset\ne A\subsetneq V(D).
\end{equation}
The strong theorem implies the weak one by the
following lemma, which we also use for finite graphs in \cref{sec:engine}.

\begin{lemma}[well-balanced implies $k$-arc-connected]\label{lem:wb-to-k}
Let $k\ge1$, let $G$ be a $2k$-edge-connected graph of any cardinality, and let $D$ be a well-balanced orientation of
$G$. Then $D$ is $k$-arc-connected; equivalently, every cut of $G$ has at least $k$ arcs of $D$ in each direction.
\end{lemma}
\begin{proof}
Let $x\ne y$. Every cut separating $x$ from $y$ has at least $2k$ edges, so \cref{lem:menger}(b) gives $2k$ edge-disjoint
$x$--$y$ paths, and $\lam_G(x,y)\ge2k$. If $\lam_G(x,y)$ is finite, \eqref{eq:wb} gives
$\lam_D(x,y)\ge\lfloor\lam_G(x,y)/2\rfloor\ge k$; if it is infinite, $\lam_D(x,y)\ge\lam_G(x,y)\ge\aleph_0$. In either
case the supremum $\lam_D(x,y)$ is at least $k$, so some family of arc-disjoint directed $x\to y$ paths has at least $k$
members: a finite supremum of a set of cardinals is its maximum, and an infinite supremum means that some family is
infinite or that there are finite families of unbounded size. Hence $D$ is $k$-arc-connected, and \eqref{eq:cutform},
applied to $A$ and to $V(G)\setminus A$, gives the statement about cuts.
\end{proof}

\section{Two properties, and why they suffice}\label{sec:roadmap}

The heart of the paper is \cref{thm:ctble}(ii), about a $2$-edge-connected graph with countably many edges. In this
section we isolate two properties of an orientation, one concerning finite cuts and one concerning infinite cuts, show
that together they make the orientation well-balanced, and explain how the rest of the paper produces an orientation
with both.

\subsection{Finite cuts}

\begin{definition}[well-balanced on finite cuts]\label{def:finwb}
An orientation $D$ of a graph $G$ is \emph{well-balanced on finite cuts} if
\[|\bd^+_D(A)|\ \ge\ \lfloor\lam_G(x,y)/2\rfloor\qquad\text{for every finite cut $\bd_G(A)$ and all $x\in A$, $y\notin A$.}\]
\end{definition}

Here $\lam_G(x,y)\le|\bd_G(A)|<\infty$, since every family of edge-disjoint $x$--$y$ paths uses distinct edges of
$\bd_G(A)$. Applying the condition to $V(G)\setminus A$ and the pair $(y,x)$ shows that also
$|\bd^-_D(A)|\ge\lfloor\lam_G(x,y)/2\rfloor$. Each instance of the condition involves only the finitely many edges of one
finite cut; this is what makes it accessible to compactness. If $G$ is $2$-edge-connected, then $\lam_G(x,y)\ge2$ for
all $x\ne y$ by \cref{lem:menger}(b), so an orientation that is well-balanced on finite cuts has at least one arc in
each direction across every finite cut.

\subsection{Infinite cuts: the demand principle}

Infinite cuts are handled by a counting principle, which is the argument of Thomassen's proof of
\cite[Thm.~11]{Th16}.

\begin{definition}[directed cycle cover; demand]\label{def:demand}
Let $D$ be a digraph. A \emph{directed cycle cover} of $D$ is a family $(C_e)_{e\in E(D)}$ in which $C_e$ is a directed
cycle of $D$ containing the arc $e$. An arc $f$ is \emph{demanded} by $e$ if $f\in E(C_e)$. The cover has \emph{finite
demand} if every arc is demanded by only finitely many arcs.
\end{definition}

\begin{theorem}[demand principle; Thomassen {\cite[proof of Thm.~11]{Th16}}]\label{thm:dp}
If a digraph $D$ of any cardinality has a directed cycle cover of finite demand, then $D$ is two-way infinite.
\end{theorem}
\begin{proof}
Let $\bd(A)$ be an infinite cut of the underlying graph and suppose that $\bd^-_D(A)$ is finite. Each
$e\in\bd^+_D(A)$ lies on the directed cycle $C_e$, which leaves $A$ along $e$ and must therefore re-enter $A$ along some
arc $f\in\bd^-_D(A)$; thus $e$ demands some $f\in\bd^-_D(A)$. Hence
$\bd^+_D(A)\subseteq\bigcup_{f\in\bd^-_D(A)}\{e:f\in E(C_e)\}$, a finite union of finite sets, so
$\bd(A)=\bd^+_D(A)\cup\bd^-_D(A)$ is finite, a contradiction. The case that $\bd^+_D(A)$ is finite is symmetric.
\end{proof}

\subsection{The cut criterion}

The two properties combine as follows.

\begin{lemma}[cut criterion]\label{lem:cuts}
Let $G$ be a graph with countably many edges, and let $D$ be an orientation of $G$ that is well-balanced on finite cuts
and two-way infinite. Then $D$ is well-balanced, and for all $x\ne y$ the supremum $\lam_D(x,y)$ is attained.
\end{lemma}
\begin{proof}
Let $x\ne y$.

\emph{Case $\lam_G(x,y)<\infty$.} Put $k:=\lfloor\lam_G(x,y)/2\rfloor$. Every $A\in\cA_D(x,y)$ satisfies
$|\bd^+_D(A)|\ge k$: by \cref{def:finwb} if $\bd_G(A)$ is finite, and because $D$ is two-way infinite if $\bd_G(A)$ is
infinite. \Cref{lem:menger}(a) gives $k$ arc-disjoint directed $x\to y$ paths, so $\lam_D(x,y)\ge k$, and the supremum
is attained by \cref{lem:lam}(a).

\emph{Case $\lam_G(x,y)=\infty$.} As $G$ has countably many edges, $\lam_G(x,y)=\aleph_0$. No finite cut separates $x$
from $y$, since a separating cut $\bd_G(A)$ bounds $\lam_G(x,y)$ by $|\bd_G(A)|$. So for every $A\in\cA_D(x,y)$ the cut
$\bd_G(A)$ is infinite and, as $D$ is two-way infinite, $\bd^+_D(A)$ is infinite. Given finitely many pairwise
arc-disjoint directed $x\to y$ paths with arc set $F$, every $A\in\cA_D(x,y)$ still has $|\bd^+_{D-F}(A)|\ge1$, so
\cref{lem:menger}(a) with $k=1$ yields a directed $x\to y$ path in $D-F$ (take the first one in a fixed enumeration of
the countably many finite paths). Starting from the empty family, this produces an infinite sequence of pairwise
arc-disjoint directed $x\to y$ paths; hence $\lam_D(x,y)=\aleph_0=\lam_G(x,y)$, and the supremum is attained.
\end{proof}

\subsection{The plan}\label{sec:plan}

By \cref{thm:dp,lem:cuts}, \cref{thm:ctble}(ii) reduces to finding an orientation $D$ of a $2$-edge-connected graph $G$
with countably many edges such that
\begin{enumerate}[label=(\Alph*)]
\item\label{plan:A} $D$ is well-balanced on finite cuts, and smooth at the vertices of finite degree, and
\item\label{plan:B} $D$ has a directed cycle cover of finite demand.
\end{enumerate}
Property~\ref{plan:A} is a conjunction of conditions each of which involves finitely many edges. Property~\ref{plan:B} is
not of this kind, and it is not preserved under limits: the set of two-way infinite
orientations is not closed (\cref{rem:notclosed}), and once lost, the property cannot be repaired by local changes
(\cref{prop:fibre}). What a limit argument can deliver is a closed set of orientations all of whose members have
property~\ref{plan:B}. The piece tree of \cref{sec:tree} converts \ref{plan:B} into finite conditions. It provides finite
$2$-edge-connected subgraphs $K(Z)$ of $G$, the \emph{cores}, and we show in \cref{sec:cycles} that \ref{plan:B} holds as
soon as
\begin{enumerate}[label=(\Alph*)]
\setcounter{enumi}{5}
\item every finite cut of $G$ contains an arc of $D$ in each direction, and
\setcounter{enumi}{2}
\item every core is strongly connected in $D$.
\end{enumerate}
By the remark after \cref{def:finwb}, (F) follows from \ref{plan:A}. So it remains to find an orientation satisfying
\ref{plan:A} and (C), and both are conjunctions of finite conditions. In \cref{sec:engine,sec:compact} they are met at
every finite level by the theorem of Kir\'aly and Szigeti, and K\H{o}nig's lemma combines the levels. \Cref{thm:ctble}(i)
then follows by a reduction to $2$-edge-connected graphs (\cref{lem:red}), and \cref{sec:all} deals with uncountable
graphs.

\section{The piece tree}\label{sec:tree}

In this section we build the structure that will confine directed cycles: a canonical rooted tree of 2EC subgraphs.
We begin with three facts about 2EC graphs, valid at any cardinality.

\subsection{Three facts about 2EC graphs}

The first fact keeps every step of the construction finite: deleting finitely many edges creates finitely many bridges.

\begin{lemma}[finitely many new bridges]\label{lem:bridges}
Let $Z$ be 2EC and let $F\subseteq E(Z)$ be finite. Then $\Br(Z-F)$ is finite.
\end{lemma}
\begin{proof}
For every pair $p,q$ of ends of edges of $F$ that lie in the same component of $Z-F$, fix one $p$--$q$ path $P_{pq}$ in
$Z-F$; there are finitely many such pairs and each path is finite. Let $b$ be a bridge of $Z-F$, let $K'$ be the
component of $Z-F$ containing $b$, and let $K_1,K_2$ be the two components of $K'-b$ (there are exactly two, one
containing each end of $b$, because $K'$ is connected and $b$ is a bridge of it). An edge of $Z$ that leaves $V(K_1)$
and is not in $F$ is an edge of $Z-F$ with an end in $K_1$, hence it is $b$, since $K_1$ is a component of $(Z-F)-b$.
So $\bd_Z(V(K_1))\subseteq\{b\}\cup F$. If $\bd_Z(V(K_1))\cap F=\emptyset$, then $\bd_Z(V(K_1))=\{b\}$ (the edge $b$
does leave $V(K_1)$), and $V(K_1)$ is a nonempty proper subset of $V(Z)$, so $b$ would be a bridge of $Z$, which is
impossible. Hence $K_1$ contains an end $p$ of an edge of $F$, and likewise $K_2$ contains an end $q$ of an edge of
$F$; both lie in $K'$. Every $p$--$q$ path of $Z-F$ lies in $K'$ and uses $b$, because $p$ and $q$ lie in different
components of $K'-b$; in particular $b\in E(P_{pq})$. Hence $\Br(Z-F)\subseteq\bigcup_{p,q}E(P_{pq})$ is finite.
\end{proof}

The second fact provides the cores.

\begin{lemma}[cores exist]\label{lem:core}
Let $Z$ be 2EC with $E(Z)\ne\emptyset$, let $T\subseteq V(Z)$ be finite and let $e\in E(Z)$. Then $Z$ has a finite 2EC
subgraph $K$ with $e\in E(K)$ and $T\subseteq V(K)$.
\end{lemma}
\begin{proof}
Let $S$ be the union of $e$ with one path in $Z$ from an end of $e$ to each $t\in T$; as $Z$ is connected these paths
exist, and $S$ is finite and connected. For each edge $g=xy$ of $S$, the graph $Z-g$ is connected because $g$ is not a
bridge of $Z$, so there is an $x$--$y$ path $P_g$ in $Z-g$, and $P_g+g$ is a cycle. Put
$K:=S\cup\bigcup_{g\in E(S)}P_g$. It is finite and connected, and every edge of $K$ lies on a cycle of $K$: an edge $g$
of $S$, and every edge of $P_g$, lies on $P_g+g\subseteq K$. So $K$ has no bridge, and $K$ is 2EC.
\end{proof}

The third fact: deleting all bridges of a graph leaves 2EC components.

\begin{lemma}[2EC components]\label{lem:classes}
For every graph $H$, every component $W$ of $H-\Br(H)$ is 2EC, and every vertex of $H$ lies in exactly one such $W$.
\end{lemma}
\begin{proof}
Let $g\in E(W)$. It is not a bridge of $H$, so it lies on a cycle $C$ of $H$. No edge of a cycle is a bridge, so
$C\subseteq H-\Br(H)$; as $C$ is connected and contains $g$, we have $C\subseteq W$, and $g$ is not a bridge of $W$.
Also $W$ is connected, and $H-\Br(H)$ has vertex set $V(H)$.
\end{proof}

We call the components of $H-\Br(H)$ the \emph{2EC components} of $H$.

\subsection{The construction}\label{sec:construction}

Throughout \cref{sec:tree,sec:cycles,sec:compact} we make the following standing assumptions: $G$ is a 2EC graph with
at least one edge and countably many edges. Since $G$ is connected and has an edge, every vertex is incident with an edge, so $V(G)$ is
countable as well. We fix an enumeration $E(G)=\{e_1,e_2,\dots\}$ (finite or infinite, without repetitions) and an
enumeration of the countably many finite subgraphs of $G$. ``The first'' object with some property always refers to
these enumerations; so all selections are canonical, and no choice principle is used in
\cref{sec:tree,sec:cycles,sec:compact}.

\begin{definition}[pieces]\label{def:piece}
A \emph{piece} is a pair $(Z,T)$ where $Z\subseteq G$ is a 2EC subgraph and $T\subseteq V(Z)$ is a finite set, whose
elements are the \emph{terminals} of the piece. We usually write $Z$ for $(Z,T)$ and $T(Z)$ for its set of terminals. If
$E(Z)\ne\emptyset$ we put $\mu(Z):=\min\{i:e_i\in E(Z)\}$, the index of the first edge of $Z$.
\end{definition}

\begin{construction}[the piece tree]\label{con:tree}
The \emph{root} is the piece $(G,\emptyset)$. A piece $(Z,T)$ with $E(Z)=\emptyset$ is a single vertex, as $Z$ is
connected; it is a \emph{leaf}, we put $K(Z):=Z$ and $B(Z):=\emptyset$, and it has no children. A piece $(Z,T)$ with
$E(Z)\ne\emptyset$ is processed as follows.
\begin{enumerate}[label=\textup{(P\arabic*)},leftmargin=3em]
\item\label{P1} \emph{Core.} $K(Z)$ is the first finite 2EC subgraph of $Z$ with $e_{\mu(Z)}\in E(K(Z))$ and
$T\subseteq V(K(Z))$; it exists by \cref{lem:core}. As $E(K(Z))\ne\emptyset$, the core has at least two vertices and
no isolated vertex.
\item\label{P2} \emph{New bridges.} $B(Z):=\Br\big(Z-E(K(Z))\big)$; this set is finite by \cref{lem:bridges}.
\item\label{P3} \emph{Children.} The children of $(Z,T)$ are the pieces $(W,T(W))$, where $W$ ranges over the
components of $\big(Z-E(K(Z))\big)-B(Z)$---each of which is 2EC by \cref{lem:classes}---and
\[T(W):=V(W)\cap\big(V(K(Z))\cup V(B(Z))\big),\]
a finite set; here $V(B(Z))$ is the set of ends of the edges in $B(Z)$.
\end{enumerate}
The pieces, defined level by level, form a rooted tree $\cT$; the \emph{depth} of a piece is its distance from the
root. We write $\cW(Z)$ for the set of children of $Z$ and
\[O(Z):=E(K(Z))\cup B(Z)\]
for the set of \emph{own edges} of $Z$ (so $O(Z)=\emptyset$ for a leaf). Since every vertex of $K(Z)$ is an end of an
edge of $K(Z)$, the set $T^*(Z):=V(K(Z))\cup V(B(Z))$ is exactly the set of ends of own edges of $Z$; it is finite.
For a non-leaf piece $Z$, the \emph{quotient} $Q(Z)$ is the finite multigraph with vertex set $\cW(Z)$ whose edges are
the own edges of $Z$ whose two ends lie in different children, each joining the two children containing its ends (own
edges with both ends in one child are discarded). An orientation of $G$ induces an orientation of $Q(Z)$.
\end{construction}

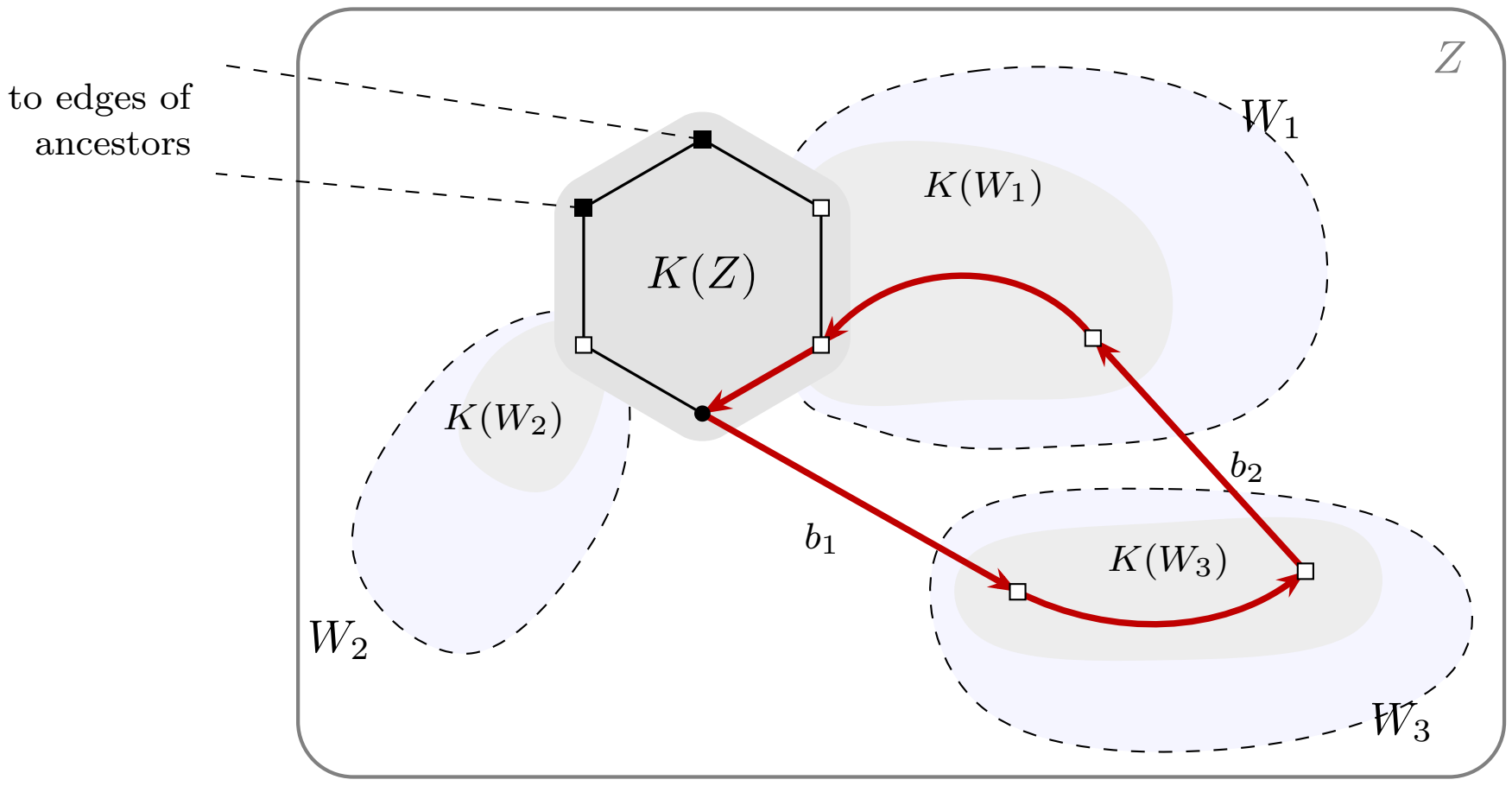
\begin{figure}[t]
\centering
\begin{tikzpicture}[x=1.05cm,y=1.05cm]
\draw[rounded corners=12pt,thick,gray] (-0.75,-3.35) rectangle (8.05,2.25);
\node[gray,font=\small] at (7.65,1.9) {$Z$};
\draw[dashed,fill=blue!4] plot[smooth cycle,tension=0.7] coordinates {(2.85,1.15) (4.2,1.8) (6.0,1.55) (6.75,0.5) (6.25,-0.65) (4.6,-0.95) (3.4,-0.8) (2.82,-0.35)};
\draw[dashed,fill=blue!4] plot[smooth cycle,tension=0.7] coordinates {(4.0,-1.55) (5.6,-1.25) (7.45,-1.6) (7.65,-2.6) (6.0,-3.15) (4.25,-2.85)};
\draw[dashed,fill=blue!4] plot[smooth cycle,tension=0.7] coordinates {(1.2,0.05) (1.65,-0.45) (1.45,-1.55) (0.5,-2.45) (-0.35,-1.7) (0.25,-0.45)};
\node[font=\small] at (6.35,1.45) {$W_1$};
\node[font=\small] at (7.3,-2.95) {$W_3$};
\node[font=\small] at (-0.45,-2.35) {$W_2$};
\fill[gray!14] plot[smooth cycle,tension=0.7] coordinates {(2.95,1.0) (4.2,1.25) (5.45,0.65) (5.45,-0.4) (4.2,-0.6) (2.95,-0.45)};
\fill[gray!14] plot[smooth cycle,tension=0.7] coordinates {(4.2,-1.7) (5.5,-1.5) (6.95,-1.55) (6.95,-2.3) (5.5,-2.5) (4.25,-2.35)};
\fill[gray!14] plot[smooth cycle,tension=0.7] coordinates {(1.25,-0.02) (1.5,-0.45) (1.1,-1.25) (0.45,-0.95) (0.65,-0.3)};
\node[font=\scriptsize] at (4.25,0.95) {$K(W_1)$};
\node[font=\scriptsize] at (5.6,-1.78) {$K(W_3)$};
\node[font=\scriptsize] at (0.75,-0.75) {$K(W_2)$};
\fill[core,rounded corners=5pt] (1.12,0.92) -- (2.2,1.55) -- (3.28,0.92) -- (3.28,-0.32) -- (2.2,-0.95) -- (1.12,-0.32) -- cycle;
\node[font=\small] at (2.2,0.3) {$K(Z)$};
\coordinate (k1) at (1.334,0.8); \coordinate (k2) at (2.2,1.3); \coordinate (k3) at (3.066,0.8);
\coordinate (k4) at (3.066,-0.2); \coordinate (k5) at (2.2,-0.7); \coordinate (k6) at (1.334,-0.2);
\draw[semithick] (k1)--(k2)--(k3)--(k4)--(k5)--(k6)--(k1);
\draw[dashed] (k1) -- (-1.35,1.05); \draw[dashed] (k2) -- (-1.35,1.85);
\node[font=\scriptsize,align=right,anchor=east] at (-1.4,1.45) {to edges of\\ancestors};
\coordinate (w) at (4.5,-2.0); \coordinate (wp) at (6.6,-1.85); \coordinate (x) at (5.05,-0.15);
\draw[cyc] (k5) -- (w) node[midway,below left,black,font=\scriptsize] {$b_1$};
\draw[cyc] (wp) -- (x) node[pos=0.45,right,black,font=\scriptsize] {$b_2$};
\draw[cyc] (w) .. controls (5.2,-2.35) and (6.0,-2.3) .. (wp);
\draw[cyc] (x) .. controls (4.6,0.45) and (3.65,0.45) .. (k4);
\draw[cyc] (k4) -- (k5);
\node[tm] at (k1) {}; \node[tm] at (k2) {};
\node[rectangle,draw,fill=white,inner sep=1.7pt] at (k3) {};
\node[rectangle,draw,fill=white,inner sep=1.7pt] at (k4) {};
\node[rectangle,draw,fill=white,inner sep=1.7pt] at (k6) {};
\node[vx] at (k5) {};
\node[rectangle,draw,fill=white,inner sep=1.7pt] at (w) {};
\node[rectangle,draw,fill=white,inner sep=1.7pt] at (wp) {};
\node[rectangle,draw,fill=white,inner sep=1.7pt] at (x) {};
\end{tikzpicture}
\caption{A schematic piece $Z$. Its core $K(Z)$ (dark) contains the terminals $T(Z)$ (black squares), where edges of
ancestors attach. Deleting $E(K(Z))$ creates the new bridges $b_1,b_2$. The children are $W_1,W_2,W_3$ (dashed) and the
single vertices on the core outside them, which are leaves. The terminals of the non-leaf children (white squares) are the
vertices where they meet $K(Z)$ or the new bridges, and they lie in the cores of the children (light). The red cycle
through $b_1$ is a local cycle as in \cref{lem:cyc}: it returns through the cores of $W_3$ and $W_1$ and through
$K(Z)$.}\label{fig:piece}
\end{figure}

No orientation is chosen in the construction. Nothing below depends on which admissible core is chosen in \ref{P1}: any
finite 2EC subgraph of $Z$ containing $e_{\mu(Z)}$ and $T$ will do, and the canonical choice only serves to avoid the
Axiom of Choice. In words: a piece keeps the finitely many edges of its core and of its new bridges, and passes every
other edge to exactly one child (\cref{lem:children}). The terminals of a child $W$ are the vertices of $W$ where own
edges of the parent attach; we shall see that they include every vertex of $W$ at which an edge owned by an ancestor
attaches (\cref{lem:bdry}). Since the core of $W$ contains its terminals, a strongly connected orientation of $K(W)$
routes directed paths between any two places where the rest of $G$ meets $W$; this is what makes the local cycles of
\cref{sec:cycles} possible. \Cref{fig:piece} shows a schematic piece.

\begin{example}[the ladder]\label{ex:ladder}
Let $\Lad$ be the ladder with vertices $a_n,b_n$ ($n\ge0$), \emph{rails} $a_na_{n+1}$, $b_nb_{n+1}$ and \emph{rungs}
$a_nb_n$; the vertices $a_n,b_n$ form \emph{level} $n$. Enumerate $E(\Lad)$ level by level, starting with $e_1=a_0b_0$.
The square on $a_0,a_1,b_1,b_0$ is the smallest admissible root core; take it. In $\Lad-E(K(\Lad))$ the vertices
$a_0,b_0$ are isolated and $a_1,b_1$ have degree one, so $B(\Lad)=\{a_1a_2,b_1b_2\}$. The children are the four
single vertices $a_0,b_0,a_1,b_1$, which are leaves, and the ladder $\Lad_{\ge2}$ on the levels $n\ge2$, with
terminals $\{a_2,b_2\}$. Inductively the core of $\Lad_{\ge2j}$ is the square on the levels $2j,2j+1$, and its new
bridges are the two rails between the levels $2j+1$ and $2j+2$ (\cref{fig:laddertree}). Every edge is a core edge or a
new bridge, and every vertex eventually becomes a leaf.
\end{example}

\begin{example}[the fan]\label{ex:fan}
Let $G$ consist of a ray $r_0r_1r_2\cdots$ and a hub $c$ joined to every $r_n$; the hub has infinite degree. With the
enumeration $cr_0,r_0r_1,cr_1,r_1r_2,cr_2,\dots$ and the triangle $cr_0r_1$ as root core, in $G-E(K(G))$ the vertex $r_0$
is isolated and $r_1$ is incident only with $r_1r_2$, which is the only new bridge. The children are the leaves $r_0,r_1$
and the fan $W_1$ on $\{c\}\cup\{r_n:n\ge2\}$, with terminals $\{c,r_2\}$. Inductively, with triangles as cores, the
core of $W_j$ (the fan on $\{c\}\cup\{r_n:n\ge2j\}$) is the triangle $cr_{2j}r_{2j+1}$ and $B(W_j)=\{r_{2j+1}r_{2j+2}\}$
(\cref{fig:fan} in \cref{sec:cycles}). The hub lies in every non-leaf piece, one of each depth, and is never a leaf,
which is harmless. The boundary
$\bd_G(V(W_j))$ consists of the $2j+1$ edges $cr_i$ ($i<2j$) and $r_{2j-1}r_{2j}$, and it meets $W_j$ exactly in its
terminals $c$ and $r_{2j}$.
\end{example}

\subsection{Basic properties}\label{sec:basic}

The next two lemmas record how edges and vertices are distributed in the tree.

\begin{lemma}[children]\label{lem:children}
Let $Z$ be a non-leaf piece.
\begin{enumerate}[label=\textup{(\alph*)}]
\item $E(Z)=O(Z)\sqcup\bigsqcup_{W\in\cW(Z)}E(W)$. The children of $Z$ are pairwise vertex-disjoint and their vertex sets
cover $V(Z)$. An edge of $Z$ whose ends lie in different children belongs to $O(Z)$.
\item $Z$ has finitely many children: $|\cW(Z)|\le\max(1,2|O(Z)|)$. Consequently there are finitely many pieces of each
depth.
\end{enumerate}
\end{lemma}
\begin{proof}
(a) By \ref{P2} and \ref{P3}, an edge of $Z$ that is not in $E(K(Z))$ is an edge of $Z-E(K(Z))$, hence either a bridge
of that graph, i.e.\ in $B(Z)$, or an edge of $(Z-E(K(Z)))-B(Z)$, which lies in exactly one component. The components
of a graph are vertex-disjoint and cover its vertex set, which here is $V(Z)$. An edge whose ends lie in different
children lies in no $E(W)$, hence in $O(Z)$.

(b) Let $\emptyset\ne\cX\subsetneq\cW(Z)$ and $U:=\bigcup_{W\in\cX}V(W)$. By (a), $U$ is a nonempty proper subset of
$V(Z)$; as $Z$ is connected, some edge of $Z$ leaves $U$, and its ends lie in different children, so it is an edge of
$Q(Z)$ joining $\cX$ to $\cW(Z)\setminus\cX$. Hence, if $|\cW(Z)|\ge2$, every child is incident with an edge of $Q(Z)$,
and $|\cW(Z)|\le2|E(Q(Z))|\le2|O(Z)|$, which is finite as $K(Z)$ is finite and $B(Z)$ is finite by
\cref{lem:bridges}. Induction on the depth gives finitely many pieces of each depth.
\end{proof}

In the construction, that $G$ has countably many edges carries structural weight only in part~(a) of the following
lemma (for its other uses see the set-theoretic remarks at the end of \cref{sec:compact}): the first edge of a piece is
always an own edge, so the index of the first edge strictly increases along every branch of the tree, and every edge is
owned at finite depth. Elsewhere in the construction the enumerations only serve to make choices canonical.

\begin{lemma}[ownership and chains]\label{lem:own}
\begin{enumerate}[label=\textup{(\alph*)}]
\item Every edge $e_m$ lies in $O(Z)$ for exactly one piece $Z=:Z_{e_m}$, its \emph{owner}, and $Z_{e_m}$ has depth
less than $m$. The pieces $Z$ with $e_m\in E(Z)$ are exactly $Z_{e_m}$ and its ancestors.
\item For every vertex $v$, the pieces containing $v$ form a chain $Z_0=G,Z_1,Z_2,\dots$ in $\cT$, with $Z_{j+1}$ a child
of $Z_j$, which is either infinite or ends with the leaf $\{v\}$. Pieces of the same depth are pairwise vertex-disjoint,
and two pieces sharing a vertex are comparable in $\cT$.
\end{enumerate}
\end{lemma}
\begin{proof}
(a) Put $Z_0:=G$. As long as $e_m\notin O(Z_j)$, the piece $Z_j$ is not a leaf (since $e_m\in E(Z_j)\ne\emptyset$), and
by \cref{lem:children}(a) there is a unique child $Z_{j+1}$ of $Z_j$ with $e_m\in E(Z_{j+1})$. Now
$E(Z_{j+1})\subseteq E(Z_j)\setminus O(Z_j)$ by \cref{lem:children}(a), and $e_{\mu(Z_j)}\in E(K(Z_j))\subseteq O(Z_j)$ by
\ref{P1}. Hence every index of an edge of $Z_{j+1}$ is at least $\mu(Z_j)$ and different from $\mu(Z_j)$, that is,
$\mu(Z_{j+1})>\mu(Z_j)$. As $e_m\in E(Z_j)$ we have $\mu(Z_j)\le m$, and $\mu(Z_0)\ge1$. So $\mu(Z_j)\ge j+1$, whence
$j+1\le m$: the sequence has at most $m$ members, and since it can only stop at a piece owning $e_m$, we get $e_m\in
O(Z_j)$ for some $j\le m-1$. Thus $Z_{e_m}:=Z_j$ has depth $j<m$.

Now let $Z$ be any piece with $e_m\in E(Z)$; we show by induction on the depth $d$ of $Z$ that $Z=Z_d$. For $d=0$, $Z$ is
the root $Z_0$. If $d\ge1$ and $P$ is the parent of $Z$, then $e_m\in E(Z)\subseteq E(P)$, so $P=Z_{d-1}$ by induction,
and $e_m\notin O(P)$ because $E(Z)\cap O(P)=\emptyset$ by \cref{lem:children}(a); so $Z$ is the unique child of
$Z_{d-1}$ containing $e_m$, i.e.\ $Z=Z_d$. Finally $e_m\in O(Z_i)$ only for $i=j$: for $i<j$ we have $e_m\in
E(Z_{i+1})$, which is disjoint from $O(Z_i)$. As a piece owning $e_m$ contains $e_m$, it is among $Z_0,\dots,Z_j$, so
the owner is unique.

(b) The root contains $v$. If a non-leaf piece contains $v$, exactly one of its children contains $v$, by
\cref{lem:children}(a). Let $Z_0=G,Z_1,\dots$ be the sequence in which $Z_{j+1}$ is the child of $Z_j$ containing $v$,
continued as long as $Z_j$ is not a leaf. Every piece containing $v$ is a member of this sequence: by induction on
its depth $d$, its parent contains $v$ and is therefore $Z_{d-1}$, and it is the child of $Z_{d-1}$ containing $v$. If
the sequence is finite, its last member is a leaf containing $v$, i.e.\ the single vertex $\{v\}$. Two pieces of the
same depth sharing a vertex are both the member of that depth of the chain of that vertex, so they coincide; and two
pieces sharing a vertex are both on the chain of that vertex, hence comparable.
\end{proof}

Next we show that every piece is attached to the rest of $G$ by finitely many edges, all owned by its ancestors, and only
at its terminals.

\begin{lemma}[boundaries and terminals]\label{lem:bdry}
\begin{enumerate}[label=\textup{(\alph*)}]
\item Every edge of $\bd_G(V(Z))$ is owned by a proper ancestor of $Z$. Hence $\bd_G(V(Z))$ is finite for every
piece~$Z$.
\item If $u\in V(Z)$ is an end of an edge owned by a proper ancestor of $Z$, then $u\in T(Z)$. If $Z$ is not a leaf, then
$T(Z)\subseteq V(K(Z))$.
\end{enumerate}
\end{lemma}
\begin{proof}
(a) Let $e=uw$ with $u\in V(Z)$ and $w\notin V(Z)$. Then $e\notin E(Z)$. The owner $Z_e$ satisfies $e\in E(Z_e)$, so
$u\in V(Z_e)$, and $Z_e$ is comparable with $Z$ by \cref{lem:own}(b). If $Z_e$ were $Z$ or a descendant of $Z$, then
$e\in E(Z_e)\subseteq E(Z)$, a contradiction; so $Z_e$ is a proper ancestor of $Z$. A piece has finitely many proper
ancestors, each owning finitely many edges.

(b) We use induction on the depth of $Z$; the root has no proper ancestor. Let $Z$ have parent $P$, and let $e\in O(A)$
be an edge at $u\in V(Z)$, where $A$ is a proper ancestor of $Z$. If $A=P$, then $u$ is an end of an edge of
$E(K(P))\cup B(P)$, so $u\in T^*(P)$; as $u\in V(Z)$, \ref{P3} gives $u\in T(Z)$. If $A$ is a proper ancestor of $P$,
then $u\in V(Z)\subseteq V(P)$, so $u\in T(P)$ by induction; $P$ is not a leaf, as it has the child $Z$, so
$T(P)\subseteq V(K(P))$ by \ref{P1}, and $u\in V(Z)\cap V(K(P))\subseteq T(Z)$ by \ref{P3}. The last sentence of~(b) is
\ref{P1}.
\end{proof}

\subsection{Finite cuts and levels}\label{sec:levels}

Only finitely many pieces are split by a finite cut, so from some depth on every piece lies on one side of it. This
allows us to read every finite cut of $G$ in a finite quotient of $G$.

\begin{lemma}[finite cuts are finite unions of pieces]\label{lem:decomp}
Let $\bd_G(X)$ be a finite cut. Call a piece $Z$ \emph{split} if $V(Z)$ meets both $X$ and $V(G)\setminus X$, and
\emph{maximal unsplit} if it is not split but its parent is.
\begin{enumerate}[label=\textup{(\alph*)}]
\item There are only finitely many split pieces; the root is split, and the parent of a split piece is split.
\item The maximal unsplit pieces are finitely many and pairwise vertex-disjoint, their vertex sets partition $V(G)$, and
$X$ is the union of the vertex sets of the maximal unsplit pieces contained in $X$.
\item If $n_0$ is the largest depth of a split piece, then for every $n\ge n_0$ no piece of depth $n+1$ is split.
\end{enumerate}
\end{lemma}
\begin{proof}
(a) A split piece $Z$ is connected and meets both sides, so it contains a path from a vertex of $X$ to a vertex outside
$X$, and some edge $f\in E(Z)$ of this path lies in $\bd_G(X)$. By \cref{lem:own}(a), $Z$ is $Z_f$ or an ancestor of
$Z_f$. There are finitely many $f\in\bd_G(X)$, each with finitely many ancestors, so there are finitely many split
pieces. The root is split because $\emptyset\ne X\ne V(G)$, and the parent of a split piece contains its vertex set,
hence is split.

(b) Let $v\in V(G)$ and consider the chain of pieces containing $v$ (\cref{lem:own}(b)). It starts with the split root,
and only finitely many of its members are split; if it is finite, its last member $\{v\}$ is unsplit. So the chain has a
first unsplit member $M(v)$, whose predecessor, its parent, is split: $M(v)$ is maximal unsplit and contains $v$. If
two maximal unsplit pieces $M\ne M'$ shared a vertex, they would be comparable, say $M$ a proper descendant of $M'$; then
the parent of $M$ is $M'$ or a descendant of $M'$, so its vertex set is contained in $V(M')$ and it is unsplit,
contradicting the maximality of $M$. Hence the maximal unsplit pieces are pairwise vertex-disjoint and cover $V(G)$.
Each of them is a child of one of the finitely many split pieces, and each piece has finitely many children
(\cref{lem:children}(b)), so there are finitely many. An unsplit piece has its vertex set inside $X$ or inside
$V(G)\setminus X$, which gives the last assertion.

(c) A piece of depth $n+1>n_0$ is not split by the choice of $n_0$.
\end{proof}

\begin{definition}[levels]\label{def:levels}
For $n\ge0$ let $E_n:=\bigcup\{O(Z):\depth(Z)\le n\}$, a finite set by \cref{lem:children}(b). Let $N_1,\dots,N_r$
($r\ge0$) be the pieces of depth $n+1$, let $L_n$ be the set of vertices of $G$ that lie in no piece of depth $n+1$, and
let
\[\cC_n:=\{V(N_1),\dots,V(N_r)\}\cup\{\{v\}:v\in L_n\}\]
be the set of \emph{classes} of level $n$. For $v\in V(G)$ we write $[v]$ for the class containing $v$ (or $[v]_n$ when
the level is not clear). $G^{(n)}$ is the finite multigraph obtained from $G$ by contracting every class to a single
vertex and deleting loops; we identify its vertex set with $\cC_n$.
\end{definition}

\begin{lemma}[properties of the levels]\label{lem:levels}
\begin{enumerate}[label=\textup{(\roman*)}]
\item $\cC_n$ is a finite partition of $V(G)$. A vertex $v$ lies in $L_n$ if and only if its chain ends with a leaf of
depth at most $n$; $L_n$ is finite, and every edge at a vertex of $L_n$ belongs to $E_n$.
\item Every edge of $G$ whose ends lie in different classes of $\cC_n$ belongs to $E_n$. Hence $E(G^{(n)})\subseteq
E_n$, and $G^{(n)}$ is finite.
\item For $\emptyset\ne\cX\subsetneq\cC_n$ and $A:=\bigcup\cX$ we have $\bd_{G^{(n)}}(\cX)=\bd_G(A)$ as edge sets; an
end of such an edge lies in $A$ if and only if its class lies in $\cX$; and $\bd_G(A)$ is a finite cut of $G$.
\item Every class of $\cC_{n+1}$ is contained in a class of $\cC_n$.
\item For every finite cut $\bd_G(X)$ there is $n_0$ such that for all $n\ge n_0$ the set $X$ is a union of classes of
$\cC_n$; then $\bd_G(X)=\bd_{G^{(n)}}(\cX_n)$ for $\cX_n:=\{C\in\cC_n:C\subseteq X\}$, with the incidences of
\textup{(iii)}.
\end{enumerate}
\end{lemma}
\begin{proof}
(i) The pieces of depth $n+1$ are finitely many and pairwise vertex-disjoint (\cref{lem:children}(b),
\cref{lem:own}(b)), so $\cC_n$ is a finite partition of $V(G)$ once $L_n$ is finite. The chain of a vertex in no piece
of depth $n+1$ does not reach depth $n+1$, hence ends with a leaf of depth at most $n$, and conversely. There are
finitely many such leaves, so $L_n$ is finite. If $v\in L_n$ and $e$ is an edge at $v$, then $Z_e$ contains $v$, so
$Z_e$ lies on the chain of $v$, whose members have depth at most $n$; hence $e\in O(Z_e)\subseteq E_n$.

(ii) Let $e$ have its ends in different classes, and suppose $e\notin E_n$. Then $Z_e$ has depth at least $n+1$, so it
has an ancestor-or-equal $N$ of depth $n+1$, and both ends of $e$ lie in $V(Z_e)\subseteq V(N)$, which is one class; a
contradiction. The edges of $G^{(n)}$ are exactly the edges of $G$ with ends in different classes.

(iii) An edge of $G$ with exactly one end in $A$ has its ends in different classes, one in $\cX$ and one outside, so it is
an edge of $G^{(n)}$ leaving $\cX$, and conversely. The statement about incidences is immediate. Finally $A$ is a
nonempty proper subset of $V(G)$, and $\bd_G(A)\subseteq E_n$ is finite by (ii).

(iv) A piece of depth $n+2$ is contained in its parent, a piece of depth $n+1$. A vertex $v\in L_{n+1}$ either lies in a
piece of depth $n+1$, whose vertex set is its class in $\cC_n$, or lies in none, and then $\{v\}$ is a class at both
levels.

(v) Let $n_0$ be as in \cref{lem:decomp}(c) and $n\ge n_0$. No piece of depth $n+1$ is split, so each class $V(N_i)$ lies
inside $X$ or inside $V(G)\setminus X$, and singleton classes do so trivially. As $X$ is a nonempty proper subset of
$V(G)$, so is $\cX_n$ of $\cC_n$, and the rest is (iii).
\end{proof}

For the ladder (\cref{ex:ladder}), $E_n$ consists of all edges with both ends on the levels $0,\dots,2n+1$ together with
the two rails between the levels $2n+1$ and $2n+2$. The pieces of depth $n+1$ are the four single vertices on the levels
$2n$ and $2n+1$ and the ladder $\Lad_{\ge2n+2}$; $L_n$ consists of the vertices on the levels below $2n$; and $G^{(n)}$ is
the finite ladder on the levels $0,\dots,2n+1$ together with one vertex for $\Lad_{\ge2n+2}$, joined to $a_{2n+1}$ and
$b_{2n+1}$.

\section{Local cycles}\label{sec:cycles}

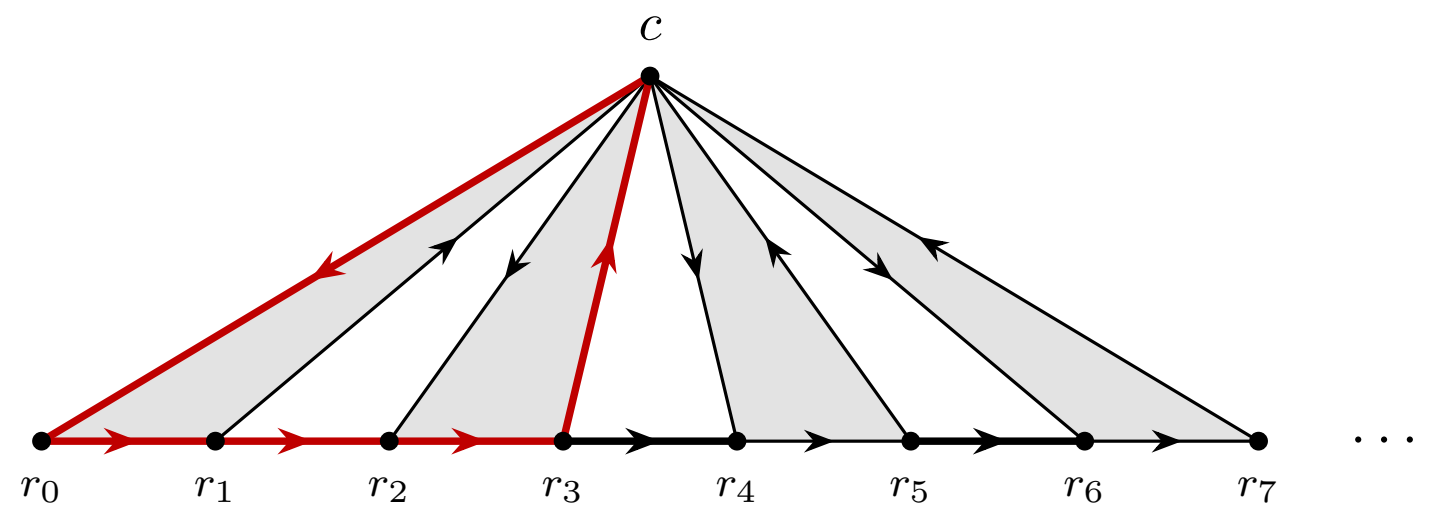
\begin{figure}[t]
\centering
\begin{tikzpicture}[x=1.15cm,y=1.15cm]
\coordinate (c) at (3.5,2.1);
\foreach \n in {0,...,7}{ \coordinate (r\n) at (\n,0); }
\foreach \j/\k in {0/1,2/3,4/5,6/7}{
  \fill[core] (c) -- (r\j) -- (r\k) -- cycle;
}
\foreach \j/\k in {4/5,6/7}{
  \draw[mid] (c) -- (r\j);
  \draw[mid] (r\j) -- (r\k);
  \draw[mid] (r\k) -- (c);
}
\draw[mid] (c) -- (r2); \draw[mid] (r1) -- (c);
\draw[midthick] (r3) -- (r4); \draw[midthick] (r5) -- (r6);
\draw[midcyc] (r1) -- (r2); \draw[midcyc] (r2) -- (r3); \draw[midcyc] (r3) -- (c); \draw[midcyc] (c) -- (r0); \draw[midcyc] (r0) -- (r1);
\foreach \n in {0,...,7}{ \node[vx] at (r\n) {}; \node[font=\scriptsize,below=3pt] at (r\n) {$r_{\n}$}; }
\node[vx] at (c) {}; \node[font=\small,above=3pt] at (c) {$c$};
\node at (7.75,0) {$\cdots$};
\end{tikzpicture}
\caption{The fan of \cref{ex:fan} with an orientation satisfying (F) and (C). The cores are the triangles
$cr_{2j}r_{2j+1}$ (shaded), the new bridges are the edges $r_{2j+1}r_{2j+2}$ (thick). The cycle of the bridge
$(r_1,r_2)$ (red) returns through the core of the child $W_1$ to the hub and through the root core back to $r_1$, inside
$\Delta(G)$. Every triangle crosses the infinite cut $\bd(\{c\})$ once in each direction.}\label{fig:fan}
\end{figure}

Throughout this section $G$ and the piece tree $\cT$ are as in \cref{sec:tree}, and $D$ is an orientation of $G$ with
the following two properties:
\begin{itemize}[leftmargin=3em]
\item[(F)] every finite cut of $G$ contains at least one arc of $D$ in each direction;
\item[(C)] $D[K(Z)]$ is strongly connected for every non-leaf piece $Z$.
\end{itemize}
We show that every arc then lies on a directed cycle that stays within the own edges of its owner and the cores of the
owner's children, and deduce that $D$ is two-way infinite (\cref{thm:demand}). \Cref{fig:laddertree,fig:piece,fig:fan}
show such cycles.

A directed cycle through a new bridge $b$ of $Z$ has to come back from the child that $b$ enters. The first step is to
see that the children of $Z$, linked by the own edges of $Z$, form a strongly connected quotient. For a union of
children that avoids the core this is where (F) is used: such a union is attached to the rest of $G$ only by new
bridges of $Z$, so its boundary is a finite cut.

\begin{lemma}[the quotient is strongly connected]\label{lem:Qstrong}
Assume \textup{(F)} and \textup{(C)}, and let $Z$ be a non-leaf piece. Let $\cK\subseteq\cW(Z)$ be the set of children
that meet $V(K(Z))$; then $\cK\ne\emptyset$.
\begin{enumerate}[label=\textup{(\alph*)}]
\item If $\emptyset\ne\cU\subseteq\cW(Z)\setminus\cK$ and $U:=\bigcup_{W\in\cU}V(W)$, then $\bd_G(U)\subseteq B(Z)$, and
$\bd_G(U)$ is exactly the set of edges of $Q(Z)$ joining $\cU$ to $\cW(Z)\setminus\cU$. Moreover, $\bd_G(U)$ is a finite
cut of $G$.
\item The orientation of $Q(Z)$ induced by $D$ is strongly connected.
\end{enumerate}
\end{lemma}
\begin{proof}
$\cK\ne\emptyset$ because $V(K(Z))\ne\emptyset$ and the children cover $V(Z)$.

(a) Let $e=xy\in\bd_G(U)$ with $x\in V(W)$ for some $W\in\cU$, and $y\notin U$. Suppose first that $e\in E(Z)$. By
\cref{lem:children}(a), $e$ lies in $E(K(Z))$, in $B(Z)$, or in $E(W')$ for a child $W'$. The first is impossible, since
$x\notin V(K(Z))$ (as $W\notin\cK$). The third is impossible: $W'=W$ would give $y\in V(W)\subseteq U$, and $W'\ne W$ would
give $x\in V(W')\cap V(W)$, contradicting the disjointness of the children. So $e\in B(Z)$. Suppose now that $e\notin
E(Z)$. The owner $Z_e$ contains $e$, so $x\in V(Z_e)$ and $Z_e$ is comparable with $Z$ (\cref{lem:own}(b)); $Z_e$ is
neither $Z$ nor a descendant of $Z$, since otherwise $e\in E(Z)$; so $Z_e$ is a proper ancestor of $Z$, and
\cref{lem:bdry}(b) gives $x\in T(Z)\subseteq V(K(Z))$, contradicting $W\notin\cK$. Hence $\bd_G(U)\subseteq B(Z)\subseteq
E(Z)$. Each edge of $\bd_G(U)$ therefore has both ends in $V(Z)$, one in a member of $\cU$ and the other in a child
outside $\cU$: it is an edge of $Q(Z)$ joining $\cU$ to its complement. Conversely, such an edge of $Q(Z)$ has exactly
one end in $U$. Finally $U\ne\emptyset$, and $U\ne V(G)$ because $U\cap V(K(Z))=\emptyset\ne V(K(Z))$; so $\bd_G(U)$ is a
cut, and it is finite since $B(Z)$ is.

(b) Let $\emptyset\ne\cY\subsetneq\cW(Z)$; we show that some arc of $Q(Z)$ leaves $\cY$ and some arc enters it. There are
three cases.

If $\cY$ contains some but not all members of $\cK$, pick $u\in V(K(Z))$ in a member of $\cY$ and $w\in V(K(Z))$ in a
member outside $\cY$. By (C) there are directed $u\to w$ and $w\to u$ paths in $D[K(Z)]$. All their vertices lie in
$V(K(Z))\subseteq V(Z)$, and hence each in exactly one child (\cref{lem:children}(a)); so the first path contains an arc of
$E(K(Z))\subseteq O(Z)$ with tail in $\bigcup\cY$ and head outside, and the second one with head in $\bigcup\cY$ and tail
outside. These arcs join different children, so they are arcs of $Q(Z)$ leaving and entering $\cY$.

If $\cY\cap\cK=\emptyset$, apply (a) to $\cU:=\cY$: $\bd_G(U)$ is a finite cut, so by (F) it contains an arc leaving $U$
and an arc entering $U$, and by (a) these are arcs of $Q(Z)$ leaving and entering $\cY$.

If $\cY\supseteq\cK$, apply (a) to $\cU:=\cW(Z)\setminus\cY\ne\emptyset$: an arc entering $\cU$ leaves $\cY$, and
conversely.

In a finite digraph in which every nonempty proper vertex set is left by an arc, every vertex reaches every other; so
$Q(Z)$ is strongly connected.
\end{proof}

For a non-leaf piece $Z$ put
\[\Delta(Z)\ :=\ O(Z)\cup\bigcup_{W\in\cW(Z)}E(K(W)),\]
the own edges of $Z$ together with the edges of the cores of its children; $\Delta(Z)$ is finite (for a leaf child $W$,
$E(K(W))=\emptyset$). The local cycles of the own edges of $Z$ live in $\Delta(Z)$.

\begin{lemma}[local cycles]\label{lem:cyc}
Assume \textup{(F)} and \textup{(C)}, and let $Z$ be a non-leaf piece.
\begin{enumerate}[label=\textup{(\alph*)}]
\item Every arc $e\in E(K(Z))$ lies on a directed cycle $C_e$ of $D$ with $E(C_e)\subseteq E(K(Z))$.
\item Every arc $b\in B(Z)$ lies on a directed cycle $C_b$ of $D$ with $E(C_b)\subseteq\Delta(Z)$.
\end{enumerate}
\end{lemma}
\begin{proof}
(a) Let $e=(u,v)$. By (C) there is a directed $v\to u$ path $P$ in $D[K(Z)]$. It does not use $e$, since $e$ enters $v$
and $v$ is the first vertex of $P$, which no arc of $P$ enters. So $P+e$ is a directed cycle inside $E(K(Z))$.

(b) The plan: the arcs of $D$ on $O(Z)$ and, as virtual arcs, all ordered pairs of distinct terminals of each child of
$Z$ form a finite digraph $\Sigma$, strongly connected by \cref{lem:Qstrong}; a directed path in $\Sigma$ from the head
of $b$ to its tail becomes a walk in $D$ when each virtual arc is replaced by a directed path in the strongly connected
core of its child. Let $T^*:=T^*(Z)=V(K(Z))\cup V(B(Z))$, the finite set of ends of own edges of $Z$. Let $\Sigma$ be the finite digraph
on $T^*$ whose arcs are the arcs of $D$ on $O(Z)$ (their ends lie in $T^*$) together with one \emph{virtual arc} $(x,y)$
for every ordered pair of distinct vertices $x,y\in T^*$ that lie in the same child $W$ of $Z$. For such a pair,
$x,y\in V(W)\cap T^*=T(W)$ by \ref{P3}; so $W$ has two distinct vertices and is not a leaf, $T(W)\subseteq V(K(W))$ by
\ref{P1}, and $D[K(W)]$ is strongly connected by (C).

\emph{$\Sigma$ is strongly connected.} Let $x,y\in T^*$ lie in the children $W_x,W_y$. By \cref{lem:Qstrong}(b) there is
a directed path in $Q(Z)$ from $W_x$ to $W_y$, that is, arcs $a_1,\dots,a_s$ of $O(Z)$ ($s\ge0$) and children
$W_x=W_0,W_1,\dots,W_s=W_y$ such that $a_i$ has its tail in $W_{i-1}$ and its head in $W_i$. The tails and heads of the
arcs $a_i\in O(Z)$ lie in $T^*$, so all virtual arcs used below exist. For $1\le i<s$ the head of $a_i$ and the tail of
$a_{i+1}$ both lie in $V(W_i)\cap T^*$; if they differ, join them by the virtual arc. Likewise join $x$ to the tail of
$a_1$, and the head of $a_s$ to $y$, by virtual arcs when they differ (if $s=0$ and $x\ne y$ use the virtual arc
$(x,y)$; if $s=0$ and $x=y$ take the trivial walk). This is a directed $x\to y$ walk in $\Sigma$.

Now let $b=(u,v)\in B(Z)$, and let $P$ be a shortest directed $v\to u$ path in $\Sigma$. As in (a), $P$ does not use
the arc $b$. Replace every virtual arc $(x,y)$ of $P$, where $x,y\in T(W)\subseteq V(K(W))$ for a child $W$, by a
directed $x\to y$ path in $D[K(W)]$. The result is a directed $v\to u$ walk $\Omega$ in $D$ with $E(\Omega)\subseteq
O(Z)\cup\bigcup_WE(K(W))=\Delta(Z)$. The walk does not use $b$: its arcs from $P$ are not $b$, and the inserted paths lie
in sets $E(K(W))\subseteq E(W)$, which do not contain $b\in B(Z)$ by \cref{lem:children}(a). Hence $u$ is reachable from
$v$ in the finite digraph formed by the arcs of $\Omega$; a shortest directed $v\to u$ path $P'$ in it, together with
$b$, is a directed cycle $C_b$ through $b$ with $E(C_b)\subseteq\Delta(Z)$.
\end{proof}

The cycles of \cref{lem:cyc} form a directed cycle cover in which every arc is demanded only by arcs owned by its own
piece or by its parent.

\begin{theorem}[finite demand]\label{thm:demand}
Assume \textup{(F)} and \textup{(C)}. For every arc $f$ of $D$ let $Z_f$ be its owner (\cref{lem:own}(a)), and let $C_f$
be the cycle given by \cref{lem:cyc}(a) if $f\in E(K(Z_f))$ and by \cref{lem:cyc}(b) if $f\in B(Z_f)$. Then:
\begin{enumerate}[label=\textup{(\roman*)}]
\item an arc $g$ is demanded only by arcs of the finite set $O(Z_g)\cup B(\parent(Z_g))$ \textup{(}the second term is
absent if $Z_g$ is the root\textup{)}; so $(C_f)_f$ is a directed cycle cover of $D$ of finite demand;
\item $D$ is two-way infinite and strongly connected;
\item if moreover every finite cut of $G$ has at least $k$ arcs of $D$ in each direction, then $D$ is $k$-arc-connected.
\end{enumerate}
\end{theorem}
\begin{proof}
(i) Let $g\in E(C_f)$. If $f\in E(K(Z_f))$, then $E(C_f)\subseteq E(K(Z_f))\subseteq O(Z_f)$, so $g\in O(Z_f)$; hence
$Z_g=Z_f$ by the uniqueness of the owner, and $f\in O(Z_g)$. If $f\in B(Z_f)$, then $E(C_f)\subseteq\Delta(Z_f)$, so
either $g\in O(Z_f)$, and again $Z_g=Z_f$ and $f\in O(Z_g)$, or $g\in E(K(W))\subseteq O(W)$ for a child $W$ of $Z_f$, in
which case $Z_g=W$, $Z_f=\parent(Z_g)$ and $f\in B(\parent(Z_g))$. Both sets are finite. Every arc has an owner, so every
arc has its cycle.

(ii) $D$ is two-way infinite by \cref{thm:dp}. Every arc lies on a directed cycle, so the two ends of every edge reach
each other in $D$; as $G$ is connected, $D$ is strongly connected.

(iii) Finite cuts have at least $k$ arcs in each direction by hypothesis, and infinite cuts have infinitely many by
(ii); conclude by \eqref{eq:cutform}.
\end{proof}

\begin{remark}[Thomassen's Theorem~11]\label{rem:th11}
Together with the orientation constructed in \cref{sec:compact}, \cref{thm:demand} gives a new proof, for graphs
with countably many edges, of Thomassen's theorem \cite[Thm.~11]{Th16} that every $2$-edge-connected graph has a strongly
connected two-way infinite orientation. Thomassen's proof builds an exhaustion $G_1\subseteq
G_2\subseteq\cdots$ by finite strongly connected oriented subgraphs and lets the cycle of an edge between $G_n$ and
$G-V(G_n)$, or of a bridge of $G-V(G_n)$, return through the accumulated $G_{n+1}$; the text asserts that every edge is
then demanded only finitely often. We could not verify this from the text for these ``attachment'' edges, since the
return path through the strongly connected $G_{n+1}$ is not constrained and may use arcs oriented at much earlier
stages. The piece tree is arranged precisely so that every return path uses only edges owned by one piece and edges of the cores of its
children, where the count of \cref{thm:demand}(i) is explicit. What is in question here is only the written
justification in \cite{Th16}, not Thomassen's procedure: no counterexample to that procedure is known to us, and we do
not claim that the theorem of \cite{Th16} is false.
\end{remark}

\section{The finite engine}\label{sec:engine}

In this section all graphs are finite; multiple edges are allowed and loops are forbidden. For such graphs Menger's
theorem gives $\lam_H(x,y)=\min\{|\bd_H(A)|:x\in A\not\ni y\}$ and $\lam_D(x,y)=\min\{|\bd^+_D(A)|:x\in A\not\ni y\}$,
and the definitions of \cref{sec:intro} specialize to those of Nash-Williams: an orientation $D$ of $H$ is
\emph{well-balanced} if $\lam_D(x,y)\ge\lfloor\lam_H(x,y)/2\rfloor$ for all $x\ne y$, \emph{smooth} if
$|\dout_D(v)-\din_D(v)|\le1$ for all $v$, and \emph{best-balanced} if it is both. As recalled in
\cref{sec:intro-finite}, every finite graph has a best-balanced orientation. We use the following generalization from \cite[p.~688]{KSz06}.

\begin{theorem}[Kir\'aly and Szigeti {\cite[Thm.~8]{KSz06}}]\label{thm:KSz}
Let $G=(V,E)$ be a graph, $\{E_1,\dots,E_k\}$ be an arbitrary partition of $E$ and let $G_i:=(V,E_i)$, $1\le i\le k$.
Then $G$ has a best-balanced orientation $\vec G$ such that the inherited orientation $\vec G_i$ of each $G_i$ is also
best-balanced.
\end{theorem}

In \cite{KSz06} graphs are finite, may have multiple edges but no loops \cite[p.~684]{KSz06}, and need not be connected
\cite[p.~685]{KSz06}; the letter $k$ in \cref{thm:KSz} only counts the parts of the partition. The theorem is also
Theorem~3.8 of \cite{BIKKS08}, where Counter-example~6 shows that the analogous statement for nested subgraphs
$H_3\subseteq H_2\subseteq H_1$ fails; we only use edge-disjoint subgraphs. The proof in \cite{KSz06} is a short argument
from Nash-Williams' odd-vertex-pairing theorem. This is where parity is dealt with; it does not appear anywhere else in
this paper. We use \cref{thm:KSz} in the following form.

\begin{lemma}[the finite statement used]\label{lem:engine}
Let $H$ be a finite graph and let $K_1,\dots,K_s$ \textup{(}$s\ge0$\textup{)} be pairwise edge-disjoint 2EC subgraphs of
$H$, each with at least one edge. Then $H$ has a best-balanced orientation $D$ in which every $D[K_j]$ is strongly
connected. If moreover $H$ is $2k$-edge-connected, then $D$ is $k$-arc-connected.
\end{lemma}
\begin{proof}
If $E(H)=\emptyset$, then $s=0$ and the empty orientation will do. Otherwise let $R:=E(H)\setminus\bigcup_jE(K_j)$ and
apply \cref{thm:KSz} to the partition $\{E(K_1),\dots,E(K_s)\}\cup\{R\}$ of $E(H)$, with $R$ omitted if it is empty. Let
$D$ be the resulting best-balanced orientation of $H$. For distinct $a,b\in V(K_j)$, the local edge-connectivity between
$a$ and $b$ in the spanning subgraph $(V(H),E(K_j))$ equals $\lam_{K_j}(a,b)$, since isolated vertices do not matter, and
$\lam_{K_j}(a,b)\ge2$, because $K_j$ is 2EC with at least two vertices, so every cut of $K_j$ has at least two edges. The
inherited orientation of $(V(H),E(K_j))$ is well-balanced, so it contains $\lfloor\lam_{K_j}(a,b)/2\rfloor\ge1$ directed
$a\to b$ paths, all inside $K_j$; hence $D[K_j]$ is strongly connected. The last sentence is \cref{lem:wb-to-k}.
\end{proof}

\begin{remark}[a polyhedral alternative for the weak theorem]\label{rem:poly}
For \cref{thm:weak} only the second conclusion of \cref{lem:engine} is needed: a $k$-arc-connected orientation of a
$2k$-edge-connected finite graph in which prescribed pairwise edge-disjoint 2EC subgraphs are strongly connected. This
can also be obtained from standard polyhedral facts without odd-vertex pairings. The in-degree vectors of the strongly
connected orientations of each $K_j$, and those of all orientations of the remaining edges, form the integer points of
integral base polyhedra that are symmetric about the respective half-degree vectors (Hakimi's theorem \cite{Ha65},
Fujishige's theorem on crossing supermodular functions and the Edmonds--Giles theorem; see \cite{Fr11}). The half-degree
vector of $H$ lies in the Minkowski sum of these polyhedra and in the polyhedron of in-degree vectors of
$k$-arc-connected orientations of $H$, and Edmonds' intersection theorem for integral base polyhedra supplies an integer
point in both, which decomposes along the parts. We do not use this route and have not written it out; we mention it
because it shows that parity is absorbed by integrality alone. The statement is also a special case of \cref{thm:KSz},
as the proof of \cref{lem:engine} shows.
\end{remark}

\section{Compactness over the levels}\label{sec:compact}

We now construct an orientation with the properties \ref{plan:A} and (C) of \cref{sec:plan}. Throughout, $G$ and $\cT$
are as in \cref{sec:tree}, with the notation of \cref{def:levels}; in addition, $V_n$ denotes the finite set of ends of
edges of $E_n$.

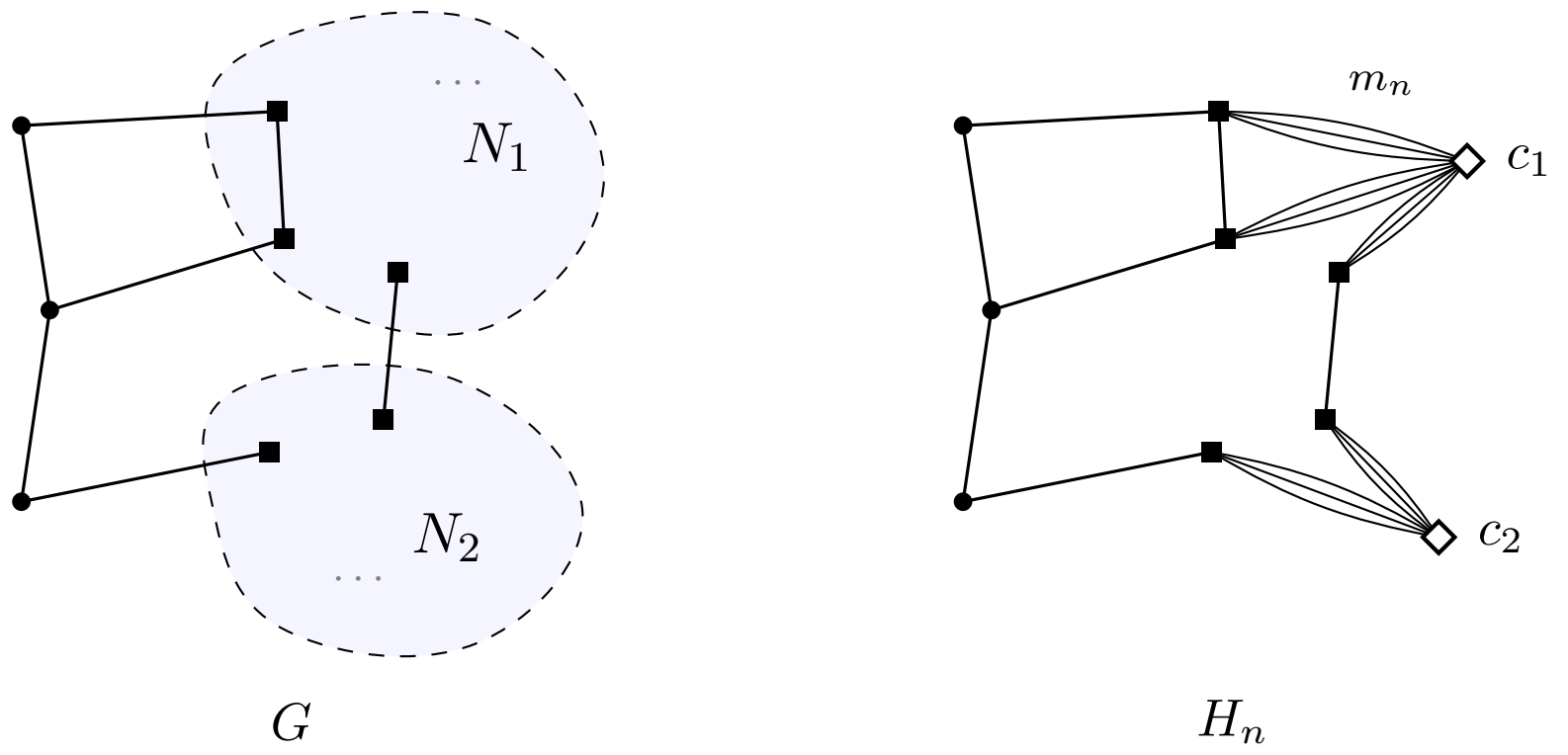
\begin{figure}[t]
\centering
\begin{tikzpicture}[x=0.95cm,y=0.95cm]
\foreach \sh/\name in {0/G,6.3/H_n}{
\begin{scope}[xshift=\sh cm]
 \coordinate (u1) at (0,1.75); \coordinate (u2) at (0.2,0.45); \coordinate (u3) at (0,-0.9);
 \coordinate (t1) at (1.8,1.85); \coordinate (t2) at (1.85,0.95); \coordinate (t3) at (2.65,0.72);
 \coordinate (t4) at (1.75,-0.55); \coordinate (t5) at (2.55,-0.32);
 \node[font=\small] at (1.9,-2.45) {$\name$};
\end{scope}}
\draw[dashed,fill=blue!4] plot[smooth cycle,tension=0.7] coordinates {(1.5,2.2) (3.2,2.5) (4.1,1.45) (3.4,0.38) (2.2,0.45) (1.45,1.2)};
\draw[dashed,fill=blue!4] plot[smooth cycle,tension=0.7] coordinates {(1.5,-0.12) (3.0,-0.02) (3.95,-0.95) (3.2,-1.9) (1.85,-1.8) (1.35,-0.95)};
\node at (3.35,1.6) {$N_1$}; \node at (3.0,-1.15) {$N_2$};
\node[font=\scriptsize,gray] at (3.1,2.05) {$\cdots$}; \node[font=\scriptsize,gray] at (2.4,-1.45) {$\cdots$};
\begin{scope}
 \coordinate (u1) at (0,1.75); \coordinate (u2) at (0.2,0.45); \coordinate (u3) at (0,-0.9);
 \coordinate (t1) at (1.8,1.85); \coordinate (t2) at (1.85,0.95); \coordinate (t3) at (2.65,0.72);
 \coordinate (t4) at (1.75,-0.55); \coordinate (t5) at (2.55,-0.32);
 \draw[semithick] (u1)--(u2)--(u3); \draw[semithick] (u1)--(t1); \draw[semithick] (u2)--(t2); \draw[semithick] (u3)--(t4);
 \draw[semithick] (t1)--(t2); \draw[semithick] (t3)--(t5);
 \foreach \u in {u1,u2,u3} \node[vx] at (\u) {};
 \foreach \t in {t1,t2,t3,t4,t5} \node[tm] at (\t) {};
\end{scope}
\begin{scope}[xshift=6.3cm]
 \coordinate (u1) at (0,1.75); \coordinate (u2) at (0.2,0.45); \coordinate (u3) at (0,-0.9);
 \coordinate (t1) at (1.8,1.85); \coordinate (t2) at (1.85,0.95); \coordinate (t3) at (2.65,0.72);
 \coordinate (t4) at (1.75,-0.55); \coordinate (t5) at (2.55,-0.32);
 \coordinate (c1) at (3.55,1.5); \coordinate (c2) at (3.35,-1.15);
 \foreach \t in {t1,t2,t3}{ \draw[thin] (c1) to[bend left=10] (\t); \draw[thin] (c1) -- (\t); \draw[thin] (c1) to[bend right=10] (\t); }
 \foreach \t in {t4,t5}{ \draw[thin] (c2) to[bend left=10] (\t); \draw[thin] (c2) -- (\t); \draw[thin] (c2) to[bend right=10] (\t); }
 \draw[semithick] (u1)--(u2)--(u3); \draw[semithick] (u1)--(t1); \draw[semithick] (u2)--(t2); \draw[semithick] (u3)--(t4);
 \draw[semithick] (t1)--(t2); \draw[semithick] (t3)--(t5);
 \foreach \u in {u1,u2,u3} \node[vx] at (\u) {};
 \foreach \t in {t1,t2,t3,t4,t5} \node[tm] at (\t) {};
 \node[gd,label={[font=\small]right:$c_1$}] at (c1) {};
 \node[gd,label={[font=\small]right:$c_2$}] at (c2) {};
 \node[font=\scriptsize] at (2.95,2.05) {$m_n$};
\end{scope}
\end{tikzpicture}
\caption{Left: the graph $G$ near level $n$. The pieces $N_1,N_2$ of depth $n+1$ (dashed) meet the edges of $E_n$
(solid) only at their terminals (squares); one edge of $E_n$ has both ends in $N_1$. The round vertices outside the
pieces are the vertices of $L_n$.
Right: the gadget graph $H_n$. It keeps $E_n$ with the true ends and replaces each $N_i$ by a gadget vertex $c_i$, joined
to every terminal of $N_i$ by $m_n$ parallel edges (three drawn). A cut of $H_n$ that separates a gadget contains $m_n$
parallel edges; a cut that separates no gadget is a finite cut of $G$ (\cref{lem:gadget}).}\label{fig:gadget}
\end{figure}

At level $n$ we orient the finite edge set $E_n$. To be able to continue the orientation at all later levels, we ask at
level $n$ for the finite-cut condition at every cut of the finite quotient $G^{(n)}$, together with strong connectivity of
the cores of depth at most $n$ and smoothness at the vertices of $L_n$ (which lie in no piece of depth $n+1$, so all
their edges lie in $E_n$). One cannot simply solve
this problem on $G^{(n)}$: contracting a piece of depth $n+1$ to a single vertex may identify distinct vertices of a
core of depth at most $n$ that lie in that piece, and may turn some core edges into loops. The \emph{gadget graph}
$H_n$ keeps all vertices of $V_n$ apart and instead attaches each piece $N_i$ of depth $n+1$ through a new vertex joined
to its terminals by many parallel edges (\cref{fig:gadget}). With enough parallel edges, cuts of $H_n$ that separate a
gadget are large, and cuts that separate no gadget are exactly the cuts of $G^{(n)}$ (\cref{lem:gadget}).

\subsection{Gadget graphs}

\begin{construction}[gadget graphs]\label{con:gadget}
For $n\ge0$ and an integer $m\ge1$, let $H_n(m)$ be the finite graph with vertex set $V_n\cup\{c_1,\dots,c_r\}$, where
$c_1,\dots,c_r$ are new vertices, one for each piece $N_i$ of depth $n+1$, and with the following edges: the edges of
$E_n$, with their ends as in $G$; and, for each $i$ and each $t\in V(N_i)\cap V_n$, $m$ parallel edges $c_it$, the
\emph{gadget edges}. The \emph{gadget} of $N_i$ is $\Gamma_i:=\{c_i\}\cup(V(N_i)\cap V_n)$. A set $Y\subseteq
V(H_n(m))$ is \emph{gadget-respecting} if for every $i$ either $\Gamma_i\subseteq Y$ or $\Gamma_i\cap Y=\emptyset$.
For $x\in V(G)$, its \emph{representative} is $\rho_n(x):=x$ if $x\in L_n$, and $\rho_n(x):=c_i$ if $x\in V(N_i)$; by
\cref{lem:gadget}(i) below, $L_n\subseteq V_n$, so $\rho_n(x)$ is a vertex of $H_n(m)$.
\end{construction}

\suppressfloats[t]

The following lemma makes this precise. Part~(iii) is the dictionary between gadget-respecting sets of $H_n(m)$ and
sets of classes of $\cC_n$; part~(iv) says that the cuts of $H_n(m)$ that respect gadgets are exactly the cuts of
$G^{(n)}$, with the same sides; and part~(v) bounds all other cuts.

\begin{lemma}[gadget identities]\label{lem:gadget}
Let $n\ge0$ and $m\ge1$.
\begin{enumerate}[label=\textup{(\roman*)}]
\item $V(N_i)\cap V_n=T(N_i)\ne\emptyset$ for each $i$, $L_n\subseteq V_n$, and $V_n=L_n\sqcup T(N_1)\sqcup\dots\sqcup
T(N_r)$. The graph $H_n(m)$ is finite and loopless.
\item For every non-leaf piece $Z$ of depth at most $n$, the core $K(Z)$ is a subgraph of $H_n(m)$, with the same
vertices, edges and incidences as in $G$. These cores are pairwise edge-disjoint 2EC subgraphs of $H_n(m)$, each with at
least one edge.
\item For a gadget-respecting set $Y$ put $\cY:=\{V(N_i):c_i\in Y\}\cup\{\{v\}:v\in L_n\cap Y\}\subseteq\cC_n$. Then
$Y\mapsto\cY$ is a bijection from the gadget-respecting sets onto the subsets of $\cC_n$; $Y=\emptyset$ if and only if
$\cY=\emptyset$, and $Y=V(H_n(m))$ if and only if $\cY=\cC_n$; a vertex $w\in V_n$ lies in $Y$ if and only if
$[w]\in\cY$; and for every $x\in V(G)$, $\rho_n(x)\in Y$ if and only if $[x]\in\cY$.
\item For every gadget-respecting set $Y$, with $\cY$ as in \textup{(iii)},
\[\bd_{H_n(m)}(Y)=\bd_{G^{(n)}}(\cY)=\bd_G\big(\textstyle\bigcup\cY\big)\qquad\text{as edge sets,}\]
and an end of such an edge lies in $Y$ if and only if it lies in $\bigcup\cY$, if and only if its class lies in $\cY$.
\item If $Y\subseteq V(H_n(m))$ is not gadget-respecting, then $\bd_{H_n(m)}(Y)$ contains $m$ parallel gadget edges.
Consequently every cut of $H_n(m)$ has at least $\min(m,\lam(G))$ edges, where $\lam(G)$ is the least size of a cut of
$G$.
\end{enumerate}
\end{lemma}
\begin{proof}
(i) Let $t\in V(N_i)$ be an end of some $e\in E_n$. The owner $Z_e$ has depth at most $n$ and contains $t$, so it is
comparable with $N_i$ (\cref{lem:own}(b)) and, having smaller depth, is a proper ancestor of $N_i$; \cref{lem:bdry}(b)
gives $t\in T(N_i)$. Conversely, let $P$ be the parent of $N_i$, a piece of depth $n$. By \ref{P3}, $T(N_i)=V(N_i)\cap
T^*(P)$ consists of ends of edges of $O(P)\subseteq E_n$, so $T(N_i)\subseteq V(N_i)\cap V_n$. Next, $T(N_i)\ne\emptyset$:
$P$ is connected and $E(P)\ne\emptyset$ (as $P$ has a child). If $V(N_i)=V(P)$, then $T(N_i)\supseteq V(K(P))\ne\emptyset$.
Otherwise $V(N_i)$ is a nonempty proper subset of $V(P)$, so some edge of $P$ leaves $V(N_i)$; by
\cref{lem:children}(a) it belongs to $O(P)$, and its end in $V(N_i)$ lies in $T^*(P)\cap V(N_i)=T(N_i)$. Next,
$L_n\subseteq V_n$: a vertex $v\in L_n$ is incident with an edge, since $G$ is connected and has at least two vertices,
and this edge lies in $E_n$ by \cref{lem:levels}(i). Every $w\in V_n$ lies in exactly one class; if $[w]=V(N_i)$ then
$w\in V(N_i)\cap V_n=T(N_i)$, and otherwise $w\in L_n$. Finally $V_n$ and $r$ are finite, and no edge of $H_n(m)$ is a
loop.

(ii) $E(K(Z))\subseteq O(Z)\subseteq E_n$, and the vertices of $K(Z)$ are ends of its edges, so $K(Z)$ is a subgraph of
$H_n(m)$ with its incidences. Two cores share no edge, since an edge has a unique owner (\cref{lem:own}(a)). Each core
is 2EC with an edge by \ref{P1}.

(iii) Let $Y$ be gadget-respecting and $w\in V_n$. If $w\in L_n$, then $w\in Y$ if and only if $\{w\}\in\cY$, by
definition. If $w\in V(N_i)$, then $w\in T(N_i)\subseteq\Gamma_i$ by (i), so $w\in Y$ iff $\Gamma_i\subseteq Y$ iff $c_i\in
Y$ iff $V(N_i)\in\cY$. The statement for $\rho_n(x)$ is the same computation. For the bijection, given
$\cX\subseteq\cC_n$ put $X_H:=\big(\bigcup\cX\cap V_n\big)\cup\{c_i:V(N_i)\in\cX\}$. If
$V(N_i)\in\cX$, then $\Gamma_i=\{c_i\}\cup T(N_i)\subseteq X_H$; if $V(N_i)\notin\cX$, then $c_i\notin X_H$ and
$T(N_i)\subseteq V(N_i)$ is disjoint from $\bigcup\cX$, so $\Gamma_i\cap X_H=\emptyset$. Thus $X_H$ is
gadget-respecting, and it is mapped to $\cX$. Conversely, a gadget-respecting $Y$ equals $X_H$ for $\cX=\cY$, by the
membership criterion just proved for $w\in V_n$ and because $c_i\in Y$ iff $V(N_i)\in\cY$. Hence the map is a bijection.
Finally $X_H=\emptyset$ for $\cX=\emptyset$, and $X_H=V_n\cup\{c_1,\dots,c_r\}=V(H_n(m))$ for $\cX=\cC_n$, using
$T(N_i)\subseteq V_n$.

(iv) Let $Y$ be gadget-respecting. A gadget edge has both ends in one $\Gamma_i$, so it never lies in $\bd_{H_n(m)}(Y)$. An edge $e\in E_n$ with ends
$w_1,w_2\in V_n$ lies in $\bd_{H_n(m)}(Y)$ iff exactly one of $w_1,w_2$ is in $Y$, iff exactly one of $[w_1],[w_2]$ is in
$\cY$ (by (iii)), iff $e\in\bd_{G^{(n)}}(\cY)$ (if $[w_1]=[w_2]$, both conditions fail). This gives the first equality and
the statement about incidences; the second equality is \cref{lem:levels}(iii) when $\emptyset\ne\cY\ne\cC_n$, and both
sides are empty otherwise.

(v) If $Y$ is not gadget-respecting, some $\Gamma_i$ meets both $Y$ and its complement; as $\Gamma_i=\{c_i\}\cup
T(N_i)$, there is $t\in T(N_i)$ with $c_i\in Y\not\ni t$ or $t\in Y\not\ni c_i$, and the $m$ parallel edges $c_it$ lie
in $\bd_{H_n(m)}(Y)$. Now let $\bd_{H_n(m)}(Y)$ be a cut, i.e.\ $\emptyset\ne Y\ne V(H_n(m))$. If $Y$ is not
gadget-respecting, the cut has at least $m$ edges; otherwise $\emptyset\ne\cY\ne\cC_n$ by (iii), and
$\bd_{H_n(m)}(Y)=\bd_G(\bigcup\cY)$ by (iv) is a cut of $G$, of size at least $\lam(G)$.
\end{proof}

\subsection{Level sets}

By \cref{lem:levels}(iii), the cuts of $G^{(n)}$ are exactly the cuts $\bd_G(A)$ of $G$ for which $A$ is a union of
classes of $\cC_n$, and these are finite. Condition~(a) below is the requirement of \cref{def:finwb} for these cuts.

\begin{definition}[pairs, multiplicities, level sets]\label{def:Sn}
Let $\cP_n$ be the set of ordered pairs $(x,y)$ of vertices of $G$ with $[x]\ne[y]$ in $\cC_n$. Put
\[m_n:=\max\Big(\{2\}\cup\big\{|\bd_G(V(N_i))|:1\le i\le r\big\}\cup\big\{\lam_G(x,y):x\ne y\in L_n\big\}\Big),\]
which is finite because $r$ and $L_n$ are finite, each $\bd_G(V(N_i))$ is finite (\cref{lem:bdry}(a)), and
$\lam_G(x,y)\le|\bd_G(\{x\})|<\infty$ for $x\ne y\in L_n$, as all edges at $x$ lie in the finite set $E_n$
(\cref{lem:levels}(i)). Let $H_n:=H_n(m_n)$. The \emph{level set} $S_n$ is the set of orientations $D_n$ of the edge set
$E_n$ (as edges of $G$, with their true ends) such that
\begin{enumerate}[label=\textup{(\alph*)}]
\item for every cut $\bd_{G^{(n)}}(\cX)$ of $G^{(n)}$, at least
$\max\{\lfloor\lam_G(x,y)/2\rfloor:x\in\bigcup\cX,\ y\notin\bigcup\cX\}$ arcs of $D_n$ \emph{leave $\cX$}, where an arc
$e$ leaves $\cX$ if $[\tail(e)]\in\cX$ and $[\head(e)]\notin\cX$;
\item $D_n[K(Z)]$ is strongly connected for every non-leaf piece $Z$ of depth at most $n$;
\item $|\dout_{D_n}(v)-\din_{D_n}(v)|\le1$ for every $v\in L_n$.
\end{enumerate}
$S_n$ is finite, since $E_n$ is.
\end{definition}

In~(a), every $\lam_G(x,y)$ with $x\in\bigcup\cX$ and $y\notin\bigcup\cX$ is at most $|\bd_G(\bigcup\cX)|$, which is
finite (\cref{lem:levels}(iii)), so the maximum exists. Moreover, $x\in\bigcup\cX$ and $y\notin\bigcup\cX$ mean exactly
that $[x]\in\cX\not\ni[y]$, and then $(x,y)\in\cP_n$. So~(a) says that for every cut $\bd_{G^{(n)}}(\cX)$ and every
$(x,y)\in\cP_n$ with $[x]\in\cX\not\ni[y]$, at least $\lfloor\lam_G(x,y)/2\rfloor$ arcs of $D_n$ leave $\cX$; in this form
it is verified and used below. The multiplicity $m_n$ is large enough for the following bound.

\begin{lemma}\label{lem:mn}
For every $(x,y)\in\cP_n$ we have $\lam_G(x,y)\le m_n$.
\end{lemma}
\begin{proof}
If $x\in V(N_i)$ for some $i$, then $y\notin V(N_i)$ (as $[x]\ne[y]$), so $V(N_i)$ contains $x$ but not $y$, and every
family of edge-disjoint $x$--$y$ paths uses distinct edges of $\bd_G(V(N_i))$; hence $\lam_G(x,y)\le|\bd_G(V(N_i))|\le
m_n$. The same holds if $y$ lies in some $N_i$. Otherwise $x,y\in L_n$, and $\lam_G(x,y)$ is one of the numbers in the
definition of $m_n$.
\end{proof}

The finite engine now shows that every level can be solved.

\begin{lemma}[$S_n\ne\emptyset$]\label{lem:nonempty}
For every $n\ge0$, the level set $S_n$ is nonempty.
\end{lemma}
\begin{proof}
\emph{Claim. For $(x,y)\in\cP_n$ we have $\rho_n(x)\ne\rho_n(y)$ and $\lam_{H_n}(\rho_n(x),\rho_n(y))\ge\lam_G(x,y)$.}
First, $\rho_n(x)\ne\rho_n(y)$ because $[x]\ne[y]$. By Menger's theorem in the finite graph $H_n$ it suffices to show
that $|\bd_{H_n}(Y)|\ge\lam_G(x,y)$ for every $Y\subseteq V(H_n)$ with $\rho_n(x)\in Y\not\ni\rho_n(y)$. If $Y$ is not
gadget-respecting, then $|\bd_{H_n}(Y)|\ge m_n\ge\lam_G(x,y)$ by \cref{lem:gadget}(v) and \cref{lem:mn}. If $Y$ is
gadget-respecting, then $[x]\in\cY\not\ni[y]$ by \cref{lem:gadget}(iii), so $A:=\bigcup\cY$ contains $x$ but not $y$, and
$\bd_{H_n}(Y)=\bd_G(A)$ by \cref{lem:gadget}(iv); every family of edge-disjoint $x$--$y$ paths in $G$ uses distinct edges of $\bd_G(A)$, so
$|\bd_{H_n}(Y)|\ge\lam_G(x,y)$. This proves the Claim.

\emph{The finite engine.} Let $Z_1,\dots,Z_s$ be the non-leaf pieces of depth at most $n$. By parts (i) and (ii)
of \cref{lem:gadget}, $H_n$ is a finite loopless graph and $K(Z_1),\dots,K(Z_s)$ are pairwise edge-disjoint 2EC
subgraphs of $H_n$, each with at least one edge. \Cref{lem:engine} gives a best-balanced orientation $D_H$ of $H_n$ in
which every $D_H[K(Z_j)]$ is strongly connected. Let $D_n$ be the restriction of $D_H$ to $E_n$. Since
$E(K(Z_j))\subseteq E_n$ and the incidences agree, $D_n[K(Z_j)]=D_H[K(Z_j)]$, so (b) holds.

\emph{Condition (a).} Let $\bd_{G^{(n)}}(\cX)$ be a cut and $(x,y)\in\cP_n$ with $[x]\in\cX\not\ni[y]$. By
\cref{lem:gadget}(iii) there is a gadget-respecting set $X_H$ with $\cY=\cX$ for $Y=X_H$; then $\rho_n(x)\in X_H\not\ni
\rho_n(y)$, and $\bd_{H_n}(X_H)=\bd_{G^{(n)}}(\cX)$ by \cref{lem:gadget}(iv), where an end of such an edge lies in $X_H$
if and only if its class lies in $\cX$. Since $D_H$ is well-balanced, the Claim gives
\[\lam_{D_H}(\rho_n(x),\rho_n(y))\ \ge\ \big\lfloor\lam_{H_n}(\rho_n(x),\rho_n(y))/2\big\rfloor\ \ge\
\lfloor\lam_G(x,y)/2\rfloor.\]
Each of the $\lam_{D_H}(\rho_n(x),\rho_n(y))$ arc-disjoint directed $\rho_n(x)\to\rho_n(y)$ paths of $D_H$ contains an
arc of $\bd^+_{D_H}(X_H)$, and distinct paths contain distinct such arcs. This is the point where the gadget
construction pays off: gadget edges never cross a gadget-respecting set, and $\bd_{H_n}(X_H)=\bd_{G^{(n)}}(\cX)\subseteq
E_n$ (\cref{lem:levels}(ii)). Hence every arc of $\bd^+_{D_H}(X_H)$ is an arc $e\in E_n$ of $D_n$ with $\tail(e)\in X_H$
and $\head(e)\notin X_H$, that is, with $[\tail(e)]\in\cX$ and $[\head(e)]\notin\cX$. So at least
$\lfloor\lam_G(x,y)/2\rfloor$ arcs of $D_n$ leave $\cX$.

\emph{Condition (c).} Let $v\in L_n$. All edges of $G$ at $v$ lie in $E_n$ (\cref{lem:levels}(i)), and $v$ is incident
with no gadget edge, since gadget edges join each $c_i$ to vertices of $V(N_i)$ and $v$ lies in no $N_i$. So the edges
of $H_n$ at $v$ are exactly the edges of $G$ at $v$, and $\dout_{D_n}(v)=\dout_{D_H}(v)$, $\din_{D_n}(v)=\din_{D_H}(v)$.
As $D_H$ is smooth, (c) holds. Thus $D_n\in S_n$.
\end{proof}

The level sets are compatible: a solution at level $n+1$ restricts to a solution at level $n$. Condition~(a)
quantifies over all cuts of $G^{(n)}$ for two reasons: every finite cut of $G$ is a cut of $G^{(n)}$ for all large $n$
(\cref{lem:levels}(v)), so that the limit orientation is well-balanced on all finite cuts; and the restriction works because every cut of $G^{(n)}$ is a cut of $G^{(n+1)}$ with the same
sides (\cref{lem:levels}(iii),(iv)).

\begin{lemma}[restriction]\label{lem:restr}
If $D_{n+1}\in S_{n+1}$, then the restriction of $D_{n+1}$ to $E_n$ lies in $S_n$.
\end{lemma}
\begin{proof}
Write $[\cdot]_n$ and $[\cdot]_{n+1}$ for the classes at the two levels, and $D_n$ for the restriction; note that
$E_n\subseteq E_{n+1}$.

(a) Let $\bd_{G^{(n)}}(\cX)$ be a cut of $G^{(n)}$ and $(x,y)\in\cP_n$ with $[x]_n\in\cX\not\ni[y]_n$. Put
$A:=\bigcup\cX$. By \cref{lem:levels}(iv), $A$ is a union of classes of $\cC_{n+1}$; let $\cX':=\{C\in\cC_{n+1}:
C\subseteq A\}$, so that $\bigcup\cX'=A$ and $\emptyset\ne\cX'\ne\cC_{n+1}$, as $A$ is a nonempty proper subset of
$V(G)$. By \cref{lem:levels}(iii) at both levels, $\bd_{G^{(n+1)}}(\cX')=\bd_G(A)=\bd_{G^{(n)}}(\cX)$ as edge sets, and
for an edge of this set an end lies in $A$ iff its class at level $n$ lies in $\cX$ iff its class at level $n+1$ lies in
$\cX'$; all these edges lie in $E_n$. Since $[x]_{n+1}\subseteq[x]_n\subseteq A$ and $[y]_{n+1}\subseteq[y]_n\subseteq
V(G)\setminus A$, we have $[x]_{n+1}\in\cX'$ and $[y]_{n+1}\notin\cX'$; in particular $[x]_{n+1}\ne[y]_{n+1}$, so
$(x,y)\in\cP_{n+1}$. Condition~(a) for $D_{n+1}$ at the cut $\bd_{G^{(n+1)}}(\cX')$ and the pair $(x,y)$ yields at
least $\lfloor\lam_G(x,y)/2\rfloor$ arcs $e$ of $D_{n+1}$ with $[\tail(e)]_{n+1}\in\cX'$ and
$[\head(e)]_{n+1}\notin\cX'$. Such an arc lies in $\bd_G(A)\subseteq E_n$, so it is an arc of $D_n$, with $\tail(e)\in A$
and $\head(e)\notin A$, i.e.\ with $[\tail(e)]_n\in\cX$ and $[\head(e)]_n\notin\cX$. So (a) holds for $D_n$.

(b) The non-leaf pieces of depth at most $n$ are among those of depth at most $n+1$, and $E(K(Z))\subseteq E_n$ for
them, so $D_n[K(Z)]=D_{n+1}[K(Z)]$ is strongly connected.

(c) A vertex in no piece of depth $n+1$ lies in no piece of depth $n+2$, so $L_n\subseteq L_{n+1}$; and all edges at
$v\in L_n$ lie in $E_n$ (\cref{lem:levels}(i)), so the in- and out-degrees of $v$ in $D_n$ and in $D_{n+1}$ agree.
\end{proof}

K\H{o}nig's lemma now combines the levels.

\begin{lemma}[K\H{o}nig's lemma]\label{lem:koenig}
Let $(S_n)_{n\ge0}$ be nonempty finite sets with maps $\pi_n:S_{n+1}\to S_n$. Then there are $s_n\in S_n$ with
$\pi_n(s_{n+1})=s_n$ for all $n$. Consequently, for the level sets of \cref{def:Sn} with the restriction maps of
\cref{lem:restr}, there is an orientation $D$ of $G$ whose restriction to $E_n$ lies in $S_n$ for every $n$.
\end{lemma}
\begin{proof}
For $m>n$ let $\pi_{m,n}:=\pi_n\circ\dots\circ\pi_{m-1}:S_m\to S_n$. Call $s\in S_n$ \emph{good} if
$s\in\pi_{m,n}(S_m)$ for every $m>n$. Good elements exist in every $S_n$: otherwise each $s\in S_n$ has some $m(s)>n$
with $s\notin\pi_{m(s),n}(S_{m(s)})$; with $M:=\max_s m(s)$ and any $t\in S_M$, the element $s:=\pi_{M,n}(t)$ satisfies
$s=\pi_{m(s),n}(\pi_{M,m(s)}(t))$, a contradiction. A good $s\in S_n$ has a good preimage in $S_{n+1}$: the preimages of
$s$ under $\pi_n$ are finitely many; if none were good, each preimage $t$ would have some $m(t)>n+1$ with
$t\notin\pi_{m(t),n+1}(S_{m(t)})$; with $M:=\max_tm(t)$, goodness of $s$ gives $u\in S_M$ with $\pi_{M,n}(u)=s$, and
$t:=\pi_{M,n+1}(u)$ is a preimage of $s$ with $t=\pi_{m(t),n+1}(\pi_{M,m(t)}(u))$, a contradiction. Now choose $s_0$ good
and inductively $s_{n+1}$ a good preimage of $s_n$, the first one in a fixed enumeration of the finite set $S_{n+1}$.
For the level sets: $E(G)=\bigcup_nE_n$ by \cref{lem:own}(a), and compatible orientations $s_n$ of the sets $E_n$ define
an orientation of $E(G)$.
\end{proof}

\subsection{Proof of \texorpdfstring{\cref{thm:ctble}}{Theorem 3}}\label{sec:proof3}

\begin{proof}[Proof of \cref{thm:ctble}\textup{(ii)}]
Let $G$ be $2$-edge-connected with countably many edges; so $G$ is 2EC with an edge, and
\cref{sec:tree,sec:cycles,sec:compact} apply. Let $D$ be the orientation given by \cref{lem:koenig}.

\emph{$D$ is well-balanced on finite cuts.} Let $\bd_G(A)$ be a finite cut, $x\in A$ and $y\notin A$; then
$\lam_G(x,y)\le|\bd_G(A)|<\infty$. By \cref{lem:levels}(v) there is $n$ such that $A$ is a union of classes of $\cC_n$
and $\bd_G(A)=\bd_{G^{(n)}}(\cX)$ for $\cX=\{C\in\cC_n:C\subseteq A\}$, an end of such an edge lying in $A$ if and only if
its class lies in $\cX$. Then $[x]\subseteq A$ and $[y]\subseteq V(G)\setminus A$, so $[x]\ne[y]$, $(x,y)\in\cP_n$, and
$[x]\in\cX\not\ni[y]$. Condition~(a) for $D|_{E_n}\in S_n$ gives at least $\lfloor\lam_G(x,y)/2\rfloor$ arcs $e$ of $D$ in
$E_n$ with $[\tail(e)]\in\cX$ and $[\head(e)]\notin\cX$, that is, with tail in $A$ and head outside $A$. Hence
$|\bd^+_D(A)|\ge\lfloor\lam_G(x,y)/2\rfloor$.

\emph{(F) and (C).} As $G$ is $2$-edge-connected, (F) follows by the remark after \cref{def:finwb}. Explicitly: let
$\bd_G(A)$ be a finite cut and pick $x\in A$, $y\notin A$; every cut separating $x$ from $y$ has at least two edges, so
$\lam_G(x,y)\ge2$, and by the previous paragraph $|\bd^+_D(A)|\ge1$; applied to the pair $(y,x)$ and the finite cut
$\bd_G(V(G)\setminus A)=\bd_G(A)$, it gives $|\bd^-_D(A)|=|\bd^+_D(V(G)\setminus A)|\ge1$. For (C), every non-leaf piece
$Z$ has some depth $n$, and $D[K(Z)]=(D|_{E_n})[K(Z)]$ is strongly connected by condition~(b) for $D|_{E_n}\in S_n$.

\emph{Two-way infinite, strongly connected, well-balanced.} By \cref{thm:demand}(ii), $D$ is two-way infinite and
strongly connected. As $D$ is also well-balanced on finite cuts, it is well-balanced by \cref{lem:cuts}. If $G$ is
$2k$-edge-connected, $D$ is $k$-arc-connected by \cref{lem:wb-to-k} (or by \cref{thm:demand}(iii)).

\emph{Smoothness.} Let $v$ have finite degree. Its finitely many edges are owned by pieces of depth less than some $d$
(\cref{lem:own}(a)). The chain of pieces containing $v$ is finite: otherwise it contains a piece $Z$ of depth $d$, which
is not a leaf as the chain continues, so $Z$ is connected with at least two vertices and some edge $e$ of $Z$ is at $v$;
but then $Z$ is $Z_e$ or an ancestor of $Z_e$ (\cref{lem:own}(a)), so $\depth(Z_e)\ge d$, a contradiction. So the chain
ends with the leaf $\{v\}$, of some depth $d'$, and $v\in L_n$ for all $n\ge d'$. Moreover all edges at $v$ lie in $E_n$
(\cref{lem:levels}(i)), so $\dout_D(v)=\dout_{D|_{E_n}}(v)$ and $\din_D(v)=\din_{D|_{E_n}}(v)$, and condition~(c) for
$D|_{E_n}\in S_n$ gives $|\dout_D(v)-\din_D(v)|\le1$. If $v$ has infinite degree, $\bd_G(\{v\})$ is an infinite cut, so it
has infinitely many arcs in each direction; as $G$ has countably many edges, $\dout_D(v)=\din_D(v)=\aleph_0$.

The second sentence of \cref{thm:ctble}(ii) holds for every smooth well-balanced two-way infinite orientation $D$ of a
$2$-edge-connected graph with countably many edges: strong connectivity because $\lam_G(x,y)\ge2$ for all $x\ne y$ and
\eqref{eq:wb}, the degrees as in the previous paragraph, and $k$-arc-connectivity by \cref{lem:wb-to-k}.
\end{proof}

To pass from $2$-edge-connected graphs to all graphs we use the 2EC components. Between different 2EC components
nothing is required, and paths between two vertices of one 2EC component stay inside it; for smoothness, we may reverse
the orientation of whole 2EC components and choose the directions of the bridges, and a compactness argument over the
tree formed by the 2EC components and the bridges does the rest.

\begin{lemma}[reduction to 2EC components]\label{lem:red}\leavevmode
\begin{enumerate}[label=\textup{(\roman*)}]
\item If every $2$-edge-connected graph with countably many edges has a well-balanced orientation, then so does every
graph with countably many edges.
\item The same holds with ``smooth well-balanced'' in place of ``well-balanced''.
\end{enumerate}
\end{lemma}
\begin{proof}
Components can be treated separately: for $x,y$ in different components $\lam_G(x,y)=0$ and nothing is required, and for
$x,y$ in one component all $x$--$y$ paths, directed or not, stay in that component. So let $G$ be connected with
countably many edges, let $\cW$ be the set of 2EC components of $G$ (each of which is 2EC by \cref{lem:classes}, has
countably many edges, and is either a single vertex or $2$-edge-connected), and for $v\in V(G)$ let $W(v)$ be the 2EC
component containing $v$.

\emph{Every edge of $G$ at a vertex of some $W\in\cW$ that is not an edge of $W$ is a bridge.} Otherwise it lies on a cycle
all of whose edges are non-bridges, hence inside one component of $G-\Br(G)$, which meets $W$ and so is $W$. In
particular every bridge joins two distinct 2EC components (a bridge with both ends in $W$ would close a cycle with a path
in the connected graph $W$), and no two bridges join the same two 2EC components (they would close a cycle with paths
inside the two components). The \emph{bridge tree}, with vertex set $\cW$ and the bridges as edges, is connected because
$G$ is, and acyclic: a cycle $b_1,\dots,b_q$ ($q\ge3$) of the bridge tree would lift, by inserting paths inside the 2EC
components, to a closed walk of $G$ that uses the edge $b_1$ exactly once, and such a walk contains a cycle through
$b_1$. So the bridge tree is a tree.

\emph{Paths between vertices of one 2EC component stay in it.} Let $x,y\in V(W)$ and let $P$ be an $x$--$y$ path of $G$
using an edge outside $E(W)$; let $b$ be the first such edge along $P$. All vertices of $P$ before $b$ lie in $V(W)$, so
$b$ has an end in $V(W)$ and is a bridge. The connected graph $W$ lies in one component of $G-b$, and after $b$ the path
is in the other component, from which it can return to $y\in V(W)$ only through $b$, which is impossible for a path.
Hence $\lam_G(x,y)=\lam_W(x,y)$ and $\lam_D(x,y)=\lam_{D[W]}(x,y)$ for every orientation $D$ of $G$.

\emph{Pairs in different 2EC components.} If $W(x)\ne W(y)$, the path of the bridge tree between them uses a bridge $b$.
The 2EC components on the two sides of $b$ in the bridge tree lie in the two components of $G-b$ (the components on one
side are joined to each other by bridges other than $b$), so $b$ separates $x$ from $y$ and $\lam_G(x,y)\le1$; the
requirement $\lfloor\lam_G(x,y)/2\rfloor=0$ is void.

(i) Orient every 2EC component with an edge well-balancedly, by hypothesis, and every bridge arbitrarily; by the two
preceding paragraphs the union is well-balanced.

(ii) Give every 2EC component $W$ with an edge a smooth well-balanced orientation $D_W$. Its reverse
$\overleftarrow{D_W}$ is also smooth and well-balanced, because $\lam_{\overleftarrow{D_W}}(x,y)=\lam_{D_W}(y,x)$ and
$\lam_W$ is symmetric. For $v\in V(W)$ of finite degree in $G$ put $\varepsilon_W(v):=\dout_{D_W}(v)-\din_{D_W}(v)\in
\{-1,0,1\}$ (for a single-vertex component put $\varepsilon_W\equiv0$). We choose a sign $s_W\in\{\pm1\}$ for every
2EC component (meaning: use $D_W$ if $s_W=1$ and its reverse if $s_W=-1$) and a direction for every bridge $b$, with
$o_b(v):=+1$ if $b$ leaves $v$ and $o_b(v):=-1$ if $b$ enters $v$, such that for every vertex $v$ of finite degree
\begin{equation}\label{eq:Sv}\tag{S$_v$}
\Big|\,s_{W(v)}\varepsilon_{W(v)}(v)+\sum_{b\in\Br(G),\,b\ni v}o_b(v)\Big|\le1 .
\end{equation}
Since every edge at $v$ is an edge of $W(v)$ or a bridge, \eqref{eq:Sv} says that the resulting orientation is smooth at
$v$, and it is well-balanced by (i). Each constraint \eqref{eq:Sv} involves finitely many of the countably many variables
$s_W$ (for the 2EC components with an edge) and $o_b$, each ranging over a two-element set. By K\H{o}nig's lemma it
suffices to show that every finite family of constraints can be satisfied: enumerate the vertices $v_1,v_2,\dots$ of
finite degree, let $S'_j$ be the finite set of assignments to the variables occurring in
$(\mathrm{S}_{v_1}),\dots,(\mathrm{S}_{v_j})$ that satisfy these constraints, with the restriction maps
$S'_{j+1}\to S'_j$, and apply \cref{lem:koenig}; variables that occur in no constraint are set arbitrarily.

So let $X$ be a finite set of vertices of finite degree. Let $\Phi$ be the graph with vertex set $\{W(x):x\in X\}$ whose
edges are the bridges with both ends in $X$, each joining the 2EC components of its two ends. Since a bridge joins two
distinct 2EC components and no two bridges join the same two, $\Phi$ is a subgraph of the bridge tree, i.e.\ a finite
forest. Call a bridge with exactly one end in $X$ \emph{dangling}; its direction influences only the constraint at its
end in $X$. Bridges with no end in $X$ are oriented arbitrarily, and $s_W$ is arbitrary for the 2EC components $W$ that
are not vertices of $\Phi$. Root each tree of $\Phi$. Every edge $b$ of $\Phi$ joins two distinct vertices of $\Phi$, one
of which is the parent of the other; so $b$ is the \emph{parent edge} $b_W$ of exactly one vertex $W$ of $\Phi$ (the
child), and a \emph{child edge} of the other at exactly one vertex of $X$, namely its end in that 2EC component. For a non-root vertex $W$ of $\Phi$ let $v_W\in
X\cap V(W)$ be the end of $b_W$ in $W$. Thus, for $x\in X\cap V(W)$, every bridge at $x$ has its other end in $X$ (an
edge of $\Phi$) or not (a dangling bridge), and \emph{the bridges at $x$ other than $b_W$ are exactly the child edges of
$\Phi$ at $x$ and the dangling bridges at $x$} (for a root $W$ there is no $b_W$). Let $c(x)$ be the number of these
bridges; it is finite since $x$ has finite degree.

We fix the variables by a top-down recursion over each tree of $\Phi$, processing a vertex $W$ of $\Phi$ only after its
parent; when a non-root $W$ is processed, the direction of $b_W$ has already been fixed (as a child edge at the parent).
\begin{enumerate}[label=\textup{(\arabic*)}]
\item At a root $W$ put $\beta_W:=0$. At a non-root $W$ put $\beta_W:=o_{b_W}(v_W)\in\{\pm1\}$, the contribution of the
already oriented parent edge at $v_W$.
\item If $W$ is not a root and $\varepsilon_W(v_W)\ne0$, choose $s_W$ with $s_W\varepsilon_W(v_W)=-\beta_W$; otherwise put
$s_W:=1$. For $x\in X\cap V(W)$ put $f(x):=s_W\varepsilon_W(x)+[x=v_W]\,\beta_W$ (the second term is absent at a root).
Then $f(x)\in\{-1,0,1\}$: for $x\ne v_W$, $f(x)=\pm\varepsilon_W(x)$; for $x=v_W$, $f(x)=0$ if $\varepsilon_W(v_W)\ne0$
and $f(x)=\beta_W$ otherwise.
\item At each $x\in X\cap V(W)$, orient the $c(x)$ child edges and dangling bridges at $x$ so that their total
contribution $\tau(x):=\sum o_b(x)$ is $0$ if $c(x)$ is even, $-f(x)$ if $c(x)$ is odd and $f(x)\ne0$, and $+1$ if $c(x)$
is odd and $f(x)=0$. This is possible because a sum of $c(x)$ terms $\pm1$ can take every value in
$\{-c(x),-c(x)+2,\dots,c(x)\}$, and each prescribed value has the parity of $c(x)$ and absolute value at most $c(x)$. In
all three cases $|f(x)+\tau(x)|\le1$. Distinct $x\in X\cap V(W)$ are treated independently: a child edge or dangling
bridge at $x$ has no other end in $V(W)$, since a bridge joins two distinct 2EC components.
\item Recurse into the children of $W$ in $\Phi$, with the directions of their parent edges just fixed in (3).
\end{enumerate}
Each vertex $W$ of $\Phi$ is processed exactly once; every edge of $\Phi$ is oriented exactly once (as a child edge
at its parent), and every dangling bridge exactly once (at the vertex of $\Phi$ containing its end in $X$). Now let $x\in X$ and
$W:=W(x)$. The left side of \eqref{eq:Sv} at $x$ is $|s_W\varepsilon_W(x)+[x=v_W]\,o_{b_W}(x)+\tau(x)|=|f(x)+\tau(x)|\le1$,
because the bridges at $x$ are $b_W$ (if $x=v_W$) together with the child edges and dangling bridges at $x$, and
$o_{b_W}(v_W)=\beta_W$. So every constraint of the finite family holds, and K\H{o}nig's lemma completes the proof.
\end{proof}

\begin{proof}[Proof of \cref{thm:ctble}\textup{(i)}]
Combine \cref{thm:ctble}(ii) with \cref{lem:red}(ii).
\end{proof}

\begin{corollary}[the weak theorem for countably many edges]\label{cor:weak-ctble}
Let $k\ge1$ and let $G$ be a $2k$-edge-connected graph with countably many edges. The orientation $D$ of
\cref{thm:ctble}\textup{(ii)} is $k$-arc-connected: every finite cut of $G$ has at least $k$ arcs of $D$ in each
direction, and every infinite cut has infinitely many.
\end{corollary}
\begin{proof}
$G$ is $2$-edge-connected, so \cref{thm:ctble}(ii) applies. The statement about cuts is \cref{lem:wb-to-k} together with
the two-way infinite property, and $k$-arc-connectivity follows from \eqref{eq:cutform}.
\end{proof}

\para{Set-theoretic remarks on \cref{thm:ctble}.} The proof uses no choice principle, provided that $G$ is given
together with an enumeration of $E(G)$ and a reference orientation, that is, an ordering of the two ends of every edge
(this is implicit in identifying orientations with points of $\{0,1\}^{E(G)}$; without it, already an enumeration of
$V(G)$ would be a countable union of two-element sets): all selections in
\cref{sec:tree,sec:cycles,sec:compact} are made along fixed enumerations, and $\lam$ is unambiguous by
\cref{lem:lam}(a). In \cref{lem:red} no orientation is chosen separately for each of the possibly infinitely many
components and 2EC components: each 2EC component $W$ with an edge inherits an enumeration of its edges from $G$, the
enumerations of its vertices and of its finite subgraphs without isolated vertices are derived from it and from the
reference orientation, and $D_W$ is the orientation that the construction of this section produces from these data; the
bridges that the proof orients arbitrarily receive their reference orientation. K\H{o}nig's lemma, in the elementary form of \cref{lem:koenig}, is used twice: for the orientation of
\cref{thm:ctble}(ii), and for smoothness across 2EC components in \cref{lem:red}(ii). That $G$ has countably many edges
carries structural weight in the construction in one place only: by \cref{lem:own}(a) every edge is owned at finite
depth, so that $E(G)=\bigcup_nE_n$ and the K\H{o}nig limit orients all of $G$. (The level sets $S_n$ are finite because
$E_n$ is, which holds at any cardinality by \cref{lem:bridges} and \cref{lem:children}(b).) Countability is also needed
in \cref{lem:cuts}, where an infinite $\lam_G$ equals $\aleph_0$, in the proof of \cref{thm:ctble}(ii) for
$\dout_D(v)=\din_D(v)=\aleph_0$, and in \cref{lem:red}(ii), where K\H{o}nig's lemma runs along an enumeration of the
vertices of finite degree. Without it \cref{lem:cuts} fails: if $x$ and $y$ are joined by $\aleph_1$ parallel edges, of
which $\aleph_1$ are oriented from $x$ to $y$ and $\aleph_0$ from $y$ to $x$, the orientation is two-way infinite and
vacuously well-balanced on finite cuts, but $\lam_D(y,x)=\aleph_0<\aleph_1=\lam_G(y,x)$. Neither local finiteness nor
ends are used:
a vertex of infinite degree lies in pieces of every depth and is served by the two-way infinite property. Infinite
parallel classes are allowed: two vertices joined by infinitely many parallel edges have $\lam_G=\aleph_0$, lie in a
common class at every level, and are served by \cref{lem:cuts}.

\section{Arbitrary cardinality}\label{sec:all}

We now remove the restriction to countably many edges. For the weak theorem we write out Thomassen's Zorn's lemma
argument \cite[\S8]{Th16} for a general constant; for the strong theorem we use the reduction theorem of Pitz and
Stegemann \cite{PS26}.

\subsection{The weak theorem: Thomassen's Zorn argument}\label{sec:zorn}

In this subsection graphs have arbitrary cardinality, and ``countable'' means countably many vertices and edges. The
argument needs two ingredients: highly edge-connected graphs contain countable highly edge-connected subgraphs through
any finite vertex set \cite[Thm.~10]{Th16}, and $k$-arc-connected orientations can be glued across a contraction.

\begin{lemma}[countable witness; cf.\ {\cite[Thm.~10]{Th16}}]\label{lem:ctble}
Let $m\in\NN$, let $G$ be an $m$-edge-connected graph of any cardinality, and let $A_0\subseteq V(G)$ be finite with
$|A_0|\ge2$. Then $G$ has a countable $m$-edge-connected subgraph $G^*$ with $A_0\subseteq V(G^*)$.
\end{lemma}
\begin{proof}
Let $G_0$ be the edgeless graph on $A_0$. Given a finite subgraph $G_j$, let $G_{j+1}$ be the union of $G_j$ with, for
each of the finitely many pairs $\{a,b\}\subseteq V(G_j)$, a set of $m$ pairwise edge-disjoint $a$--$b$ paths of $G$;
these exist by \cref{lem:menger}(b), since every cut of $G$ has at least $m$ edges. Each $G_{j+1}$ is finite (it consists
of finitely many finite paths, whatever the sizes of the parallel classes), so $G^*:=\bigcup_jG_j$ is countable. Let
$\emptyset\ne A\subsetneq V(G^*)$ and pick $x\in A$ and $y\in V(G^*)\setminus A$. Both lie in some $G_j$, and the $m$
edge-disjoint $x$--$y$ paths placed in $G_{j+1}\subseteq G^*$ each contain an edge of $\bd_{G^*}(A)$. So $G^*$ is
$m$-edge-connected. (The recursion makes countably many dependent selections of path systems; this uses the axiom of
Dependent Choice.)
\end{proof}

\begin{lemma}[gluing across a contraction]\label{lem:glue}
Let $k\ge1$, let $G_1\subseteq G$ be a subgraph with a $k$-arc-connected orientation $D_1$, and let $X:=V(G_1)$. Let
$H:=G/X$ be the graph obtained from $G$ by deleting all edges with both ends in $X$ and identifying $X$ to a new vertex
$v_0$; so $E(H)\subseteq E(G)$. Let $G_2\subseteq H$ be a subgraph with $v_0\in V(G_2)$, and let $D_2$ be a
$k$-arc-connected orientation of $G_2$. Read each arc of $D_2$ at $v_0$ as the arc of $G$ at its original end in $X$.
Then $D:=D_1\cup D_2$ is a $k$-arc-connected orientation of the subgraph $G_1\cup G_2$ of $G$, whose vertex set is
$X\cup(V(G_2)\setminus\{v_0\})$.
\end{lemma}
\begin{proof}
$E(G_1)\cap E(G_2)=\emptyset$, since the edges of $G_1$ have both ends in $X$ and are not edges of $H$; so $D$ is a
well-defined orientation. We verify \eqref{eq:cutform}. Let $\emptyset\ne A\subsetneq V(D)$.

(i) If $\emptyset\ne A\cap X\ne X$: an arc of $D_1$ leaving $A\cap X$ has its tail in $A$ and its head in $X\setminus
A$, so $\bd^+_{D_1}(A\cap X)\subseteq\bd^+_D(A)$, and the left side has at least $k$ arcs.

(ii) If $X\subseteq A$: put $A':=(A\setminus X)\cup\{v_0\}\subseteq V(G_2)$. Since $A\ne V(D)$, some vertex of
$V(G_2)\setminus\{v_0\}$ lies outside $A$, so $\emptyset\ne A'\subsetneq V(G_2)$. An arc of $D_2$ leaving $A'$ has its
head in $V(G_2)\setminus A'\subseteq V(D)\setminus A$, and its tail in $A'$, which read in $G$ is a vertex of $A$ (a tail
$v_0$ becomes a vertex of $X\subseteq A$). So $|\bd^+_D(A)|\ge|\bd^+_{D_2}(A')|\ge k$.

(iii) If $A\cap X=\emptyset$: then $A\subseteq V(G_2)\setminus\{v_0\}$, so $\emptyset\ne A\subsetneq V(G_2)$, and an arc of
$D_2$ leaving $A$ has its head outside $A$ also when read in $G$ (a head $v_0$ becomes a vertex of $X$, which is disjoint
from $A$). So $|\bd^+_D(A)|\ge|\bd^+_{D_2}(A)|\ge k$.
\end{proof}

\begin{theorem}[from countable to arbitrary graphs; Thomassen {\cite[\S8]{Th16}}]\label{thm:U}
Let $k\ge1$. If every countable $2k$-edge-connected graph has a $k$-arc-connected orientation, then every
$2k$-edge-connected graph of any cardinality has one.
\end{theorem}
\begin{proof}
Let $G$ be $2k$-edge-connected; so $|V(G)|\ge2$. By \cref{lem:ctble}, with $m=2k$ and any two vertices as $A_0$, $G$ has
a countable $2k$-edge-connected subgraph $G_0$, and by hypothesis $G_0$ has a $k$-arc-connected orientation $D_0$. Let
$\cP$ be the set of pairs $(G',D')$ with $G_0\subseteq G'\subseteq G$ and $D'$ a $k$-arc-connected orientation of $G'$
extending $D_0$, ordered by $(G',D')\le(G'',D'')$ if $G'\subseteq G''$ and $D''$ extends $D'$. $\cP$ contains
$(G_0,D_0)$. For a nonempty chain $\cC\subseteq\cP$, the union $(\bigcup_{\cC}G',\bigcup_{\cC}D')$ is an oriented
subgraph of $G$ extending $D_0$ (the orientations agree on common edges), and it is $k$-arc-connected: any two of its
vertices lie in a common member $(G',D')$ of the chain, which contains $k$ arc-disjoint directed paths between them, and
these persist in the union. So every chain has an upper bound, and by \emph{Zorn's lemma} $\cP$ has a maximal element
$(G_1,D_1)$.

\emph{$G_1$ is an induced subgraph.} If an edge $e=uv\in E(G)\setminus E(G_1)$ had $u,v\in V(G_1)$, adding $e$ to $D_1$
in either direction would give a larger member of $\cP$, as adding arcs destroys no directed path.

\emph{$V(G_1)=V(G)$.} Otherwise let $X:=V(G_1)$ and $H:=G/X$ as in \cref{lem:glue}. $H$ is loopless by construction, as
every edge with both ends in $X$ is deleted. $H$ is $2k$-edge-connected: $|V(H)|\ge2$ since $X\ne V(G)$, and for
$\emptyset\ne A\subsetneq V(H)$ let $\hat A$ be $A$ with $v_0$ replaced by the set $X$ if $v_0\in A$; then
$\emptyset\ne\hat A\subsetneq V(G)$ and $\bd_H(A)=\bd_G(\hat A)$ as edge sets (an edge of $G$ with both ends in $X$ lies on
one side of $\hat A$ and is not an edge of $H$; every other edge of $G$ is an edge of $H$ and crosses $A$ if and only if
it crosses $\hat A$), so $|\bd_H(A)|\ge2k$. Apply \cref{lem:ctble} to $H$ with $A_0=\{v_0,w\}$ for some $w\in
V(H)\setminus\{v_0\}$: it gives a countable $2k$-edge-connected subgraph $G_2\subseteq H$ containing $v_0$ and $w$. (The
parallel classes of $H$ at $v_0$ may be uncountable even if $G$ is simple; this does not matter.) By hypothesis $G_2$ has
a $k$-arc-connected orientation $D_2$, and \cref{lem:glue} shows that $(G_1\cup G_2,D_1\cup D_2)\in\cP$. It is strictly
larger than $(G_1,D_1)$ because $w\notin X$, contradicting maximality.

So $G_1$ is induced and spanning, i.e.\ $G_1=G$, and $D_1$ is the required orientation.
\end{proof}

\begin{proof}[Proof of \cref{thm:weak}]
A countable $2k$-edge-connected graph has countably many edges, so it has a $k$-arc-connected orientation by
\cref{cor:weak-ctble}. \Cref{thm:U} extends this to all cardinalities.
\end{proof}

\para{Set-theoretic remarks on \cref{thm:weak}.} Besides the choice-free proof of the countable case, the proof uses
Dependent Choice in \cref{lem:ctble} and Zorn's lemma in \cref{thm:U}. It uses neither \cref{lem:lam}(b) nor
bond-faithful decompositions. The two-way infinite property is not transferred to uncountable graphs
(\cref{rem:two-routes}).

\subsection{The strong theorem: the reduction of Pitz and Stegemann}\label{sec:PS}

\begin{theorem}[Pitz and Stegemann {\cite[Thm.~1.2]{PS26}}]\label{thm:PS}
Let $\cA$ be a class of graphs that is closed under subgraphs. If every countable graph in $\cA$ has a well-balanced
orientation, then all graphs in $\cA$ have a well-balanced orientation.
\end{theorem}

In \cite{PS26}, graphs may have parallel edges but no loops (\S2 there); ``countable'' means countably many vertices and
edges (footnote~6 there, in the statement of Theorem~6.2); and well-balanced is Definition~2.1 there, that is,
\eqref{eq:wb}. For the countable graphs to which the hypothesis is applied, \cref{lem:lam}(a) shows that \eqref{eq:wb} is
unambiguous without any choice principle; for the graphs of arbitrary cardinality in the conclusion, \cref{lem:lam}(b)
shows that every $\lam_G(x,y)$ and $\lam_D(x,y)$ is attained and equals the least size of a separating cut, so that
well-balancedness in the sense of \eqref{eq:wb} coincides with \cite[Def.~2.1]{PS26} at every cardinality. The proof of
\cref{thm:PS} (\S6 of \cite{PS26}) uses the theorem of Laviolette \cite{La05} in the form proved by Thomassen
\cite[Thm.~3]{Th17}: every graph has an edge-decomposition into connected countable subgraphs $H$ such that every finite
minimal cut of $H$ is a cut of $G$. It uses the Axiom of Choice. We have not re-derived the proof of \cref{thm:PS} and use
it as a black box.

\begin{proof}[Proof of \cref{thm:strong}]
Let $\cA$ be the class of all graphs. A countable graph has countably many edges, and so has a well-balanced
orientation by \cref{thm:ctble}(i); by \cref{lem:lam}(a) this is a well-balanced orientation in the sense of
\cite[Def.~2.1]{PS26}. \Cref{thm:PS} gives a well-balanced orientation in the sense of \cite[Def.~2.1]{PS26} of every
graph, and by \cref{lem:lam}(b) this is a well-balanced orientation in the sense of \eqref{eq:wb}.
\end{proof}

\para{Set-theoretic remarks on \cref{thm:strong}.} The proof uses the Axiom of Choice, through \cref{thm:PS} and
\cite{Th17}, and through \cref{lem:lam}(b), which identifies the readings of \eqref{eq:wb} for uncountable graphs.

\begin{remark}[two routes to the weak theorem]\label{rem:two-routes}
\Cref{thm:weak} also follows from \cref{thm:strong} by \cref{lem:wb-to-k}. The route through \cref{thm:U} is
self-contained (Zorn's lemma, Dependent Choice and the countable case), whereas the route through \cref{thm:PS} depends on
\cite{PS26} and \cite{Th17}. The two-way infinite property is not transferred to uncountable graphs by either route: it
is not preserved under unions of chains, and distinct parts of a bond-faithful decomposition may share vertices.
Likewise smoothness (\cref{thm:ctble}(i)) is not claimed beyond countably many edges; see \cref{rem:smooth-open}.
\end{remark}

\section{Examples and remarks}\label{sec:examples}

\subsection{The two-way infinite property must be built in}

\begin{remark}[two-way infinite is not a closed condition]\label{rem:notclosed}
The set of two-way infinite orientations is not closed in the product topology on $\{0,1\}^{E(G)}$, so it cannot be
imposed by compactness directly. This is why the closed conditions (F) and (C), each a conjunction of clopen conditions, are imposed instead, and the two-way
infinite property is recovered by \cref{thm:demand}. Example: on the ladder $\Lad$ of \cref{ex:ladder} let $\Dlad$ orient
$a_n\to a_{n+1}$, $b_{n+1}\to b_n$, $b_0\to a_0$, and $a_n\to b_n$ for $n\ge1$ (\cref{fig:ladder}). It is strongly
connected: from any vertex one reaches $b_0$ (along the bottom rail, which a top vertex $a_m$ reaches through a rung
$a_n\to b_n$ with $n\ge\max(m,1)$), and from $b_0$ one reaches $a_0$, every top vertex, and through the rungs every bottom
vertex. But the infinite cut $\bd(\{a_n:n\ge0\})$, which consists of all rungs, has only one arc entering
$\{a_n:n\ge0\}$. For $N\ge1$ let $D^{(N)}$ agree with $\Dlad$ on all rails and on the rungs at the levels $\le N$, and
orient the rung at level $n>N$ as $b_n\to a_n$ if $n-N$ is odd and as $a_n\to b_n$ if $n-N$ is even. Then
$D^{(N)}\to\Dlad$ coordinatewise as $N\to\infty$, and each $D^{(N)}$ is two-way infinite. Indeed, every arc with both
ends at levels $>N$ lies on a directed cycle within four consecutive levels: a rung $a_n\to b_n$ (so $n\ge N+2$) on the
square $a_{n-1}a_nb_nb_{n-1}$; a rung $b_n\to a_n$ on the square $a_na_{n+1}b_{n+1}b_n$; and the rails $a_na_{n+1}$ and
$b_{n+1}b_n$ on the latter square if $n-N$ is odd, and on the octagon
$a_na_{n+1}a_{n+2}b_{n+2}b_{n+1}b_nb_{n-1}a_{n-1}$ if $n-N$ is even. The finitely many remaining arcs lie on directed
cycles because $D^{(N)}$ is strongly connected (by the argument given for $\Dlad$, except that a bottom vertex $b_n$
with $n>N$ and $n-N$ odd is reached not through its own rung but via the downward rung $a_{n+1}\to b_{n+1}$ and then the
bottom rail), and the resulting cycle cover has finite demand, so \cref{thm:dp} applies.
\end{remark}

In the $k$-fold ladder, the deficiency of $\Dlad$ cannot be repaired afterwards by any change that alters the out-degrees
at only finitely many vertices.

\begin{proposition}[the two-way infinite property cannot be repaired]\label{prop:fibre}
Let $k\ge1$, let $k\Lad$ be the ladder with every edge replaced by $k$ parallel edges, and let $D_k$ be the orientation of
$k\Lad$ that orients all $k$ copies of each edge as $\Dlad$ does. Then $k\Lad$ is $2k$-edge-connected, $D_k$ is
$k$-arc-connected, and every orientation $D'$ of $k\Lad$ with $\dout_{D'}(v)=\dout_{D_k}(v)$ for all but finitely many
vertices $v$ has only finitely many rungs oriented from $b$ to $a$; hence $D'$ is not two-way infinite. In particular,
reversing the arcs of any family of pairwise arc-disjoint directed cycles and double rays, and changing the orientation
of finitely many edges, never turns $D_k$ into a two-way infinite orientation.
\end{proposition}
\begin{proof}
$\Lad$ is connected and every edge lies on a square, so every cut of $\Lad$ has at least two edges and every cut of $k\Lad$
at least $2k$. $\Dlad$ is strongly connected (\cref{rem:notclosed}), so every cut has an arc of $\Dlad$ in each direction
and at least $k$ arcs of $D_k$ in each direction; thus $D_k$ is $k$-arc-connected by \eqref{eq:cutform}. Let $U$ be the
finite set of vertices $v$ with $\dout_{D'}(v)\ne\dout_{D_k}(v)$, and $c:=\sum_{v\in U}|\dout_{D'}(v)-\dout_{D_k}(v)|$. For
a finite set $A$ of vertices of a locally finite graph and any orientation, summing $\dout(v)-\din(v)$ over $A$ counts
the arcs inside $A$ zero times, so $|\bd^+(A)|-|\bd^-(A)|=\sum_{v\in A}(2\dout(v)-d(v))$, that is,
$|\bd^+(A)|=\frac12\big(|\bd(A)|+\sum_{v\in A}(2\dout(v)-d(v))\big)$. Apply this to $Y_m:=\{b_0,\dots,b_m\}$ for both
orientations: $|\bd^+_{D'}(Y_m)|\le|\bd^+_{D_k}(Y_m)|+c$. In $D_k$ the only arcs leaving $Y_m$ are the $k$ copies of
$b_0\to a_0$ (the rungs at the levels $1,\dots,m$ and the rails $b_{m+1}b_m$ point into $Y_m$), so
$|\bd^+_{D'}(Y_m)|\le k+c$. Every rung at a level $\le m$ that $D'$ orients from $b$ to $a$ leaves $Y_m$. As $m$ is
arbitrary, $D'$ orients at most $k+c$ rungs from $b$ to $a$, and the infinite cut $\bd(\{a_n:n\ge0\})$, which consists of
all rungs, has at most $k+c$ arcs entering $\{a_n:n\ge0\}$. For the last sentence: reversing a subdigraph in which every
vertex has equal in- and out-degree does not change any out-degree, and changing finitely many edges changes the
out-degrees at finitely many vertices only.
\end{proof}

In \cref{thm:ctble}, by contrast, the two-way infinite property is present from the start, through the cores. For
$k\Lad$ the requirement \eqref{eq:wb} is, for instance, $k$ between consecutive vertices of a rail and
$\lfloor3k/2\rfloor$ between the two ends of a rung at a level $n\ge1$; it is delivered by condition~(a) at the finite
cuts of the quotients $G^{(n)}$.

\subsection{Test cases}

\begin{example}[the ladder, continued]\label{ex:ladder2}
In $\Lad$ every cut has at least two edges and every vertex has degree at most three, so $\lam(x,y)\in\{2,3\}$ for all
$x\ne y$, and a well-balanced orientation of $\Lad$ is the same as a strongly connected one; \cref{thm:ctble}(ii) gives a
strongly connected two-way infinite orientation, which is the statement of \cite[Thm.~11]{Th16} for $\Lad$. In an
orientation with (F) and (C), the new bridges $a_1a_2$, $b_1b_2$ are oriented in opposite directions, by (F) at the
finite cut $\{a_1a_2,b_1b_2\}$, and the cycle of the bridge $a_1a_2$ given by \cref{lem:cyc} runs through the core of
$\Lad_{\ge2}$, back along $b_2b_1$ and through the root core, inside $\Delta(\Lad)$ (\cref{fig:laddertree}). An arc of the
root core is demanded only by own edges of the root, a set that does not depend on the depth of the tree.
\end{example}

\begin{example}[the twin binary tree]\label{ex:twin}
Let $T$ be the rooted binary tree, whose vertices are the finite $0$--$1$ sequences; take two copies $T\times\{0\}$ and
$T\times\{1\}$, and join $(t,0)$ to $(t,1)$ and $(t^\frown0,0)$ to $(t^\frown1,1)$ for every $t$ \cite[Fig.~1]{AKP25}.
This locally finite graph has uncountably many ends, and it is the example showing that the decomposition theorem
\cite[Thm.~2.2]{AKP25} needs countably many ends. It is $2$-edge-connected (every edge lies on a cycle of length four),
and for the piece tree it is simply a graph with countably many edges: there are finitely many pieces of each depth,
each with finitely many terminals; nothing in \cref{sec:tree,sec:cycles,sec:compact} refers to its ends; and K\H{o}nig's
lemma is applied to the finite level sets $S_n$, not to the end space.
\end{example}

\begin{example}[Aharoni--Thomassen graphs]\label{ex:AT}
For every $k$ there are countable graphs of finite edge-connectivity at least $k$ without two edge-disjoint spanning trees
\cite{AT89}; by the construction there they are locally finite (cf.\ \cite{AKP25}). The same examples show that in
infinite digraphs $k$-arc-connectivity does not yield $k$ arc-disjoint spanning arborescences. The proof here never uses spanning
trees, arborescences or packings: the cores are finite 2EC subgraphs, found by \cref{lem:core} from local paths only.
\end{example}

\begin{example}[the countable complete graph]\label{ex:Kaleph}
In $K_{\aleph_0}$ there are no finite cuts, $\lam(x,y)=\aleph_0$ for every pair, and \eqref{eq:wb} demands
$\lam_D(x,y)=\aleph_0$. Whatever finite 2EC root core $K$ is chosen (a triangle, say), $K_{\aleph_0}-E(K)$ is still
connected and bridgeless, since any two vertices have infinitely many common neighbours; so $B=\emptyset$, every piece has
exactly one child, whose vertex set is all of $V$ and whose terminals are the vertices of the previous core. Hence
$V(K(Z_0))\subseteq V(K(Z_1))\subseteq\cdots$, and every vertex lies in all but finitely many of the pairwise
edge-disjoint strongly connected cores: $\dout_D(v)=\din_D(v)=\aleph_0$ directly, and every infinite cut separates two
vertices that lie in infinitely many common cores, each of which crosses the cut in both directions. Here $G^{(n)}$ is a
single vertex, $\cP_n=\emptyset$, $H_n$ consists of $E_n$ and one gadget, and the content of the proof is (C) together
with \cref{thm:dp}.
\end{example}

\begin{example}[infinite parallel classes]\label{ex:parallel}
Let $u$ and $v$ be joined by $\aleph_0$ parallel edges, and attach a ladder at $v$. Then $\lam_G(u,v)=\aleph_0$, the pair
lies in no $\cP_n$, and the two-way infinite property serves it (\cref{lem:cuts}). If $e_1$ is a $uv$-edge and the root
core consists of two parallel $uv$-edges, then $B=\emptyset$ and the only child has terminals $\{u,v\}$. The chain of $u$
is infinite, since every piece of it contains infinitely many $uv$-edges and owns only finitely many of them, and all
$uv$-edges that a piece does not own lie in its child containing $u$ (\cref{lem:children}(a)); and $u$
lies in the core of every piece of its chain: it lies in the root core, and if it lies in the core of a piece of its
chain, then it is a terminal of the next piece of its chain by \ref{P3}, hence lies in its core by \ref{P1}. These cores are pairwise edge-disjoint and strongly connected, so $\bd(\{u\})$ has infinitely many
arcs in each direction. An uncountable parallel class is outside \cref{thm:ctble} and is handled by
\cref{thm:weak,thm:strong}: in the proof of \cref{thm:U} the countable witness uses countably many of the parallel edges.
\end{example}

\subsection{Further remarks}

\begin{remark}[comparison with earlier approaches]\label{rem:buys}
(1) \emph{No decompositions.} The proofs in \cite{Th16,As25,AKP25,AR25} decompose $V(G)$ into a finite set and finitely
many boundary-linked sets with finite boundaries, whose boundary edges start pairwise edge-disjoint equivalent rays, and
orient an increasing sequence of finite immersions; the version in which these sets are the components of $G-A$ requires
hypotheses on the ends \cite[Thm.~2.2, Fig.~1]{AKP25}. The piece tree is a decomposition of a different kind: the pieces
are nested, infinitely many of them may be infinite, and nothing links their boundaries. Infinite cuts are handled
afterwards, by the demand count, which only needs that every arc has a return cycle of bounded ``age'' in the tree.
(2) \emph{No ends.} Nothing in \cref{sec:tree,sec:cycles,sec:compact} mentions ends, rays or local finiteness; a vertex
of infinite degree simply lies in pieces of every depth. (3) \emph{No parity.} At connectivity exactly $2k$, lifting
arguments meet the parity of odd boundaries (the exceptional configuration of the lifting graph, \cite[Thm.~2.3]{As25}).
We never lift or split off, and parity enters only inside \cref{thm:KSz}, a finite theorem applied to one finite graph
at each level. (4) \emph{Why strongly connected 2EC cores.} The finite constraints must keep the finite problems solvable
at every level. Strong connectivity of pairwise edge-disjoint 2EC cores is exactly what \cref{thm:KSz} delivers for free,
and a core must contain prescribed vertices, its terminals, which a $2$-edge-connected graph can always accommodate
(\cref{lem:core}). Prescribing instead a fixed consistent orientation of a closed Eulerian subgraph through the
terminals (the closed case of \cite[Thm.~5.1]{AKP25}, which also allows open Eulerian subgraphs) is not always possible:
$K_{2,3}$ has no connected subgraph with all degrees even that contains all five vertices.
\end{remark}

\begin{remark}[smoothness for uncountable graphs]\label{rem:smooth-open}
\Cref{thm:ctble}(i) gives smooth well-balanced orientations of graphs with countably many edges. Whether every graph has
a smooth well-balanced orientation we leave open: the reduction of \cref{thm:PS} orients the parts of a bond-faithful
decomposition separately, and distinct parts may share a vertex of finite degree, where the imbalances add up; the
argument of \cref{lem:red}(ii) used that 2EC components are vertex-disjoint. The two-way infinite property is likewise not
claimed beyond countably many edges (\cref{rem:two-routes}).
\end{remark}

\begin{remark}[a transfer principle; heuristic]\label{rem:transfer}
The infinite part of the argument (\cref{sec:tree,sec:cycles,sec:compact}) uses about the orientation property only that
it is (i) decided by finitely many finite cuts at a time, (ii) inherited by the quotients $G^{(n)}$ and by restriction
from level $n+1$ to level $n$, and (iii) realizable in finite graphs simultaneously with strong connectivity of
prescribed pairwise edge-disjoint 2EC subgraphs, provided the gadget multiplicity is chosen so that the gadget graph does
not lower any finite local connectivity that a cut of $G^{(n)}$ can see (\cref{lem:mn} and the Claim in
\cref{lem:nonempty}). It seems likely that every requirement with these three features holds for all graphs with countably
many edges; the two orientation theorems of Nash-Williams are the two instances treated here. We have not formalized this
statement and do not use it.
\end{remark}

\appendix
\section{Menger's theorem in infinite graphs}\label{app:menger}

In this appendix we prove \cref{lem:menger,lem:lam}. The proofs are the standard augmenting path arguments; no
compactness is used, and the Axiom of Choice is used only in \cref{lem:lam}(b).

\begin{proof}[Proof of \cref{lem:menger}\textup{(a)}]
Necessity: each directed $x\to y$ path uses an arc of $\bd^+_D(A)$, and distinct paths of an arc-disjoint family use
distinct arcs. For sufficiency we show by induction on $j<k$ that from $j$ arc-disjoint directed $x\to y$ paths
$P_1,\dots,P_j$ we can obtain $j+1$ such paths (not necessarily containing the $P_i$). Let $F$ be the finite set of arcs
of $P_1,\dots,P_j$. The \emph{residual digraph} $R$ has vertex set $V(D)$ and consists of all arcs of $D$ not in $F$
together with a formal reversal $\bar a$ of every $a\in F$. Let $A$ be the set of vertices reachable from $x$ by a
directed walk in $R$.

\emph{Case $y\notin A$.} Let $a=(p,q)$ be an arc of $D$ with $p\in A$ and $q\notin A$. If $a\notin F$, then $a\in R$ and
$q$ would be reachable; so $a\in F$. Next let $a=(q,p)\in F$ with $q\notin A$ and $p\in A$; then $\bar a=(p,q)\in R$, and
again $q$ would be reachable. So every arc of $D$ leaving $A$ belongs to $F$, and no arc of $F$ enters $A$. Each $P_i$
starts at $x\in A$ and ends at $y\notin A$, so it contains an arc leaving $A$; after that arc it never returns to $A$,
since this would require an arc of $F$ entering $A$. Hence each $P_i$ contains exactly one arc of $\bd^+_D(A)$, every arc
of $\bd^+_D(A)$ lies on some $P_i$, and $|\bd^+_D(A)|=j<k$, contradicting the hypothesis.

\emph{Case $y\in A$.} Let $Q$ be a directed $x\to y$ path in $R$. Regard $F$ as a $0/1$-valued function on the arcs of
$D$; its net out-flow (arcs out minus arcs in) is $j$ at $x$, $-j$ at $y$ and $0$ at every other vertex, since each $P_i$
contributes $+1$ at $x$, $-1$ at $y$ and $0$ elsewhere. Let $F'$ be obtained from $F$ by adding every arc of $Q$ that is
an arc of $D$ (these are not in $F$) and deleting every $a\in F$ with $\bar a\in Q$. Adding an arc $(p,q)$ changes the net
out-flow by $+1$ at $p$ and $-1$ at $q$, and deleting $a=(q,p)$ has the same effect. So traversing $Q$ changes the net
out-flow by $+1$ at $x$, by $-1$ at $y$ and by $0$ elsewhere, and $F'$ is a $0/1$-function with net out-flow $j+1$ at
$x$, $-(j+1)$ at $y$ and $0$ elsewhere. Let $D'$ be the finite digraph whose arcs are the arcs of $F'$ and whose vertices
are their ends together with $x$ and $y$. For $A'\subseteq V(D')$ with $x\in A'$ and $y\notin A'$, summing the net
out-flow over $A'$ gives $|\bd^+_{D'}(A')|-|\bd^-_{D'}(A')|=j+1$, so $|\bd^+_{D'}(A')|\ge j+1$. By Menger's theorem for
finite digraphs (arc version; see \cite[Ch.~9]{Sch03}), $D'\subseteq D$ contains $j+1$ arc-disjoint directed $x\to y$
paths.
\end{proof}

\begin{proof}[Proof of \cref{lem:menger}\textup{(b)}]
Let $\vec G$ be the digraph obtained by replacing each edge $e=pq$ by two arcs $e^+=(p,q)$ and $e^-=(q,p)$. Then
$|\bd^+_{\vec G}(A)|=|\bd_G(A)|\ge k$ for all $A$ as in the statement, so part~(a) gives $k$ arc-disjoint directed $x\to y$
paths in $\vec G$. Among all such families choose one, $P_1,\dots,P_k$, of least total length (total lengths are natural
numbers). Suppose some edge $e=pq$ had $e^+$ on $P_i$ and $e^-$ on $P_j$. Then $i\ne j$, since a path traversing
$p\to q\to p$ would repeat $p$. Write $P_i=P_i[x,p]\,e^+\,P_i[q,y]$ and $P_j=P_j[x,q]\,e^-\,P_j[p,y]$. The directed
$x\to y$ walks $P_i[x,p]P_j[p,y]$ and $P_j[x,q]P_i[q,y]$ use disjoint arc sets, both contained in $(E(P_i)\cup
E(P_j))\setminus\{e^+,e^-\}$; shortcutting each to a path and replacing $P_i,P_j$ by them gives a family of $k$
arc-disjoint directed $x\to y$ paths of smaller total length, a contradiction. Hence no edge is used in both directions,
and the underlying undirected paths of $P_1,\dots,P_k$ are edge-disjoint.
\end{proof}

\begin{proof}[Proof of \cref{lem:lam}\textup{(a)}]
As $D$ has countably many arcs, every $|\bd^+_D(A)|$ lies in $\NN\cup\{\aleph_0\}$, so the minimum $\kappa$ exists. No
family has more than $\kappa$ members, since distinct members use distinct arcs of $\bd^+_D(A)$ for every
$A\in\cA_D(x,y)$. If $\kappa<\infty$, \cref{lem:menger}(a) provides $\kappa$ paths. If $\kappa=\aleph_0$, then
$\bd^+_D(A)$ is infinite for every $A\in\cA_D(x,y)$. Given finitely many arc-disjoint directed $x\to y$ paths with arc
set $F$, the digraph $D-F$ still has $|\bd^+_{D-F}(A)|=\infty\ge1$ for all such $A$, so \cref{lem:menger}(a) with $k=1$
yields a directed $x\to y$ path in $D-F$; we take the first one in a fixed enumeration of the countably many finite
paths, so that no choice is made. Starting from the empty family, this produces an infinite sequence of pairwise
arc-disjoint directed $x\to y$ paths. The statement for graphs follows in the same way from \cref{lem:menger}(b).
\end{proof}

\begin{proof}[Proof of \cref{lem:lam}\textup{(b)}]
The set of cardinals is nonempty, since $\{x\}\in\cA_D(x,y)$, and it has a least element because the cardinals are well
ordered (this uses the Axiom of Choice). Fix $A_0\in\cA_D(x,y)$ with $|\bd^+_D(A_0)|=\kappa$.

\emph{No family has more than $\kappa$ members.} Every directed $x\to y$ path starts in $A_0$ and ends outside $A_0$, so
it contains an arc of $\bd^+_D(A_0)$. Assigning to each member of a family $\cP$ of pairwise arc-disjoint directed
$x\to y$ paths the first arc of $\bd^+_D(A_0)$ on it defines a map $\cP\to\bd^+_D(A_0)$, which is injective because
distinct members share no arc. Hence $|\cP|\le\kappa$.

\emph{A family of $\kappa$ members exists.} If $\kappa$ is finite, \cref{lem:menger}(a) with $k:=\kappa$ applies, since
$|\bd^+_D(A)|\ge\kappa$ for every $A\in\cA_D(x,y)$ by the choice of $\kappa$ (for $\kappa=0$ the empty family will do).
Let $\kappa$ be infinite, identified with its initial ordinal, so that every ordinal $\alpha<\kappa$ has $|\alpha|<\kappa$.
Fix a well-ordering $\preceq$ of the set of all directed $x\to y$ paths of $D$ (Axiom of Choice). We define $P_\alpha$ for
$\alpha<\kappa$ by transfinite recursion. Suppose that pairwise arc-disjoint directed $x\to y$ paths $P_\beta$
($\beta<\alpha$) have been defined, and let $F_\alpha:=\bigcup_{\beta<\alpha}E(P_\beta)$. Each $E(P_\beta)$ is finite, so
$F_\alpha$ is finite if $\alpha$ is finite, and $|F_\alpha|\le|\alpha|\cdot\aleph_0=|\alpha|$ if $\alpha$ is infinite; in
both cases $|F_\alpha|<\kappa$, as $\kappa$ is infinite and $|\alpha|<\kappa$. (The cruder bound $|\alpha|\cdot\aleph_0$ is
not smaller than $\kappa$ when $\kappa=\aleph_0$ and $\alpha\ge1$, so the finite case must be treated separately.) Now let
$A\in\cA_D(x,y)$. If $\bd^+_D(A)\subseteq F_\alpha$, then $\kappa\le|\bd^+_D(A)|\le|F_\alpha|<\kappa$, which is absurd;
so some arc of $\bd^+_D(A)$ is not in $F_\alpha$, i.e.\ $|\bd^+_{D-F_\alpha}(A)|\ge1$, where $D-F_\alpha$ has vertex set
$V(D)$ and arc set $E(D)\setminus F_\alpha$. By \cref{lem:menger}(a) with $k=1$, applied to $D-F_\alpha$, there is a
directed $x\to y$ path in $D-F_\alpha$; let $P_\alpha$ be the $\preceq$-least one. It is arc-disjoint from every
$P_\beta$, $\beta<\alpha$, so the recursion continues. The resulting paths $P_\alpha$ ($\alpha<\kappa$) are pairwise
arc-disjoint, and pairwise distinct because each has at least one arc ($x\ne y$). So $\{P_\alpha:\alpha<\kappa\}$ is a
family of exactly $\kappa$ pairwise arc-disjoint directed $x\to y$ paths.

Together: the cardinalities of families of pairwise arc-disjoint directed $x\to y$ paths are bounded by $\kappa$, and
$\kappa$ occurs among them, so their supremum $\lam_D(x,y)$ equals $\kappa$ and is attained. For a graph $G$ the same proof
applies verbatim, with $\bd_G$ in place of $\bd^+_D$ (every $x$--$y$ path contains an edge of $\bd_G(A_0)$), the family of
sets $A\subseteq V(G)$ with $x\in A$ and $y\notin A$ in place of $\cA_D(x,y)$, and \cref{lem:menger}(b) (with $k:=\kappa$,
and with $k=1$ applied to $G-F_\alpha$) in place of \cref{lem:menger}(a).
\end{proof}

\section*{Acknowledgement}

The authors used the AI system Claude (Anthropic) extensively in this work, both to develop the proofs and to write the
manuscript, and take full responsibility for its content.



\begin{thebibliography}{BIKKS08}
\bibitem[AT89]{AT89} R.~Aharoni and C.~Thomassen, \emph{Infinite, highly connected digraphs with no two arc-disjoint
spanning trees}, J.~Graph Theory \textbf{13} (1989), 71--74.
\bibitem[As22]{As22} A.~Assem, \emph{Analysis of the lifting graph}, arXiv:2212.03347 (2022; v3, 2024).
\bibitem[As25]{As25} A.~Assem, \emph{Towards Nash-Williams orientation conjecture for infinite graphs}, J.~Graph Theory
\textbf{108} (2025), no.~3, 608--619; doi:10.1002/jgt.23192; arXiv:2306.10631.
\bibitem[AKP25]{AKP25} A.~Assem, M.~Koloschin, and M.~Pitz, \emph{The Nash-Williams orientation theorem for graphs with
countably many ends}, European J.~Combin.\ \textbf{124} (2025), 104043; doi:10.1016/j.ejc.2024.104043;
arXiv:2310.03601.
\bibitem[AR25]{AR25} A.~Assem and R.~B.~Richter, \emph{Edge-disjoint linkage in infinite graphs}, J.~Combin.\ Theory
Ser.~B \textbf{182} (2027), 173--193; doi:10.1016/j.jctb.2026.09.002.
\bibitem[AMP26]{AMP26} L.~Aurichi, P.~Magalh\~aes Jr., and G.~E.~Pinto, \emph{On orientations preserving
edge-connectivity in infinite graphs}, J.~Logic Comput.\ \textbf{36} (2026), no.~5, exag028; doi:10.1093/logcom/exag028;
arXiv:2510.06449.
\bibitem[BK10]{BK10} J.~Bar\'at and M.~Kriesell, \emph{What is on his mind?}, Discrete Math.\ \textbf{310} (2010),
2573--2583.
\bibitem[BIKKS08]{BIKKS08} A.~Bern\'ath, S.~Iwata, T.~Kir\'aly, Z.~Kir\'aly, and Z.~Szigeti, \emph{Recent results on
well-balanced orientations}, Discrete Optim.\ \textbf{5} (2008), 663--676; doi:10.1016/j.disopt.2008.03.001.
(Numbering of statements as in the technical report version, EGRES TR-2006-11.)
\bibitem[Egy41]{Egy41} L.~Egyed, \emph{\"Uber die wohlgerichteten unendlichen Graphen}, Mat.\ Fiz.\ Lapok \textbf{48}
(1941), 505--509 (in Hungarian, with German summary).
\bibitem[Fr11]{Fr11} A.~Frank, \emph{Connections in Combinatorial Optimization}, Oxford Lecture Series in Mathematics
and its Applications~38, Oxford University Press, Oxford, 2011.
\bibitem[Ha65]{Ha65} S.~L.~Hakimi, \emph{On the degrees of the vertices of a directed graph}, J.~Franklin Inst.\
\textbf{279} (1965), 290--308.
\bibitem[Hu91]{Hu91} A.~Huck, \emph{A sufficient condition for graphs to be weakly $k$-linked}, Graphs Combin.\
\textbf{7} (1991), 323--351.
\bibitem[Jo03]{Jo03} T.~Jord\'an, \emph{Constrained edge-splitting problems}, SIAM J.~Discrete Math.\ \textbf{17}
(2003), 88--102.
\bibitem[KSz06]{KSz06} Z.~Kir\'aly and Z.~Szigeti, \emph{Simultaneous well-balanced orientations of graphs},
J.~Combin.\ Theory Ser.~B \textbf{96} (2006), 684--692; doi:10.1016/j.jctb.2006.01.002.
\bibitem[KP20]{KP20} M.~Koloschin and M.~Pitz, \emph{Orientations of infinite graphs}, arXiv:2011.05936 (2020).
\bibitem[La05]{La05} F.~Laviolette, \emph{Decompositions of infinite graphs: I---bond-faithful decompositions},
J.~Combin.\ Theory Ser.~B \textbf{94} (2005), 259--277.
\bibitem[NW60]{NW60} C.~St.~J.~A.~Nash-Williams, \emph{On orientations, connectivity and odd-vertex-pairings in finite
graphs}, Canad.\ J.~Math.\ \textbf{12} (1960), 555--567; doi:10.4153/CJM-1960-049-6.
\bibitem[NW67]{NW67} C.~St.~J.~A.~Nash-Williams, \emph{Infinite graphs---a survey}, J.~Combin.\ Theory \textbf{3}
(1967), 286--301.
\bibitem[NW69]{NW69} C.~St.~J.~A.~Nash-Williams, \emph{Well-balanced orientations of finite graphs and unobtrusive
odd-vertex-pairings}, in: W.~T.~Tutte (ed.), Recent Progress in Combinatorics (Proc.\ Third Waterloo Conf.\ on
Combinatorics, May 1968), Academic Press, New York, 1969, pp.~133--149.
\bibitem[ORT16]{ORT16} S.~Ok, R.~B.~Richter, and C.~Thomassen, \emph{Liftings in finite graphs and linkages in infinite
graphs with prescribed edge-connectivity}, Graphs Combin.\ \textbf{32} (2016), 2575--2589.
\bibitem[PS26]{PS26} M.~Pitz and J.~Stegemann, \emph{The strong Nash-Williams orientation theorem for rayless graphs},
J.~Graph Theory (2026), published online 23 July 2026, doi:10.1002/jgt.70100.
\bibitem[Rob39]{Rob39} H.~E.~Robbins, \emph{A theorem on graphs, with an application to a problem of traffic control},
Amer.\ Math.\ Monthly \textbf{46} (1939), 281--283.
\bibitem[Sch03]{Sch03} A.~Schrijver, \emph{Combinatorial Optimization: Polyhedra and Efficiency}, Algorithms and
Combinatorics~24, Springer, Berlin, 2003.
\bibitem[Th16]{Th16} C.~Thomassen, \emph{Orientations of infinite graphs with prescribed edge-connectivity},
Combinatorica \textbf{36} (2016), 601--621; doi:10.1007/s00493-015-3173-0.
\bibitem[Th17]{Th17} C.~Thomassen, \emph{Nash-Williams' cycle-decomposition theorem}, Combinatorica \textbf{37} (2017),
1027--1037; doi:10.1007/s00493-016-3424-8.
\end{thebibliography}
\end{document}